\documentclass{ws-m3as}

\begin{document}

\markboth{Y. Jia, Y. Du and L. Guo}
{Geometry-dependent Ekman layer approximations on curved domains}

%
\catchline{}{}{}{}{}
%

\title{Non-flat Ekman Boundary Layers: Topographic Lift, Generalized Ekman Pumping, and Anisotropic Asymptotic Behavior
}

\author{Yifei Jia}
\address{College of Mathematics and System Sciences,\\
Xinjiang University,\\
Urumqi 830017, P. R. China\\
jiayifei333@163.com}

\author{Yi Du\footnote{Corresponding author}}
\address{Department of Mathematics, Jinan University,\\
Guangzhou 510632, P. R. China\\
duyidy@jnu.edu.cn}

\author{Lihui Guo}
\address{College of Mathematics and System Sciences,\\
	Xinjiang University,\\
	Urumqi 830017, P. R. China\\
	lihguo@126.com}

\maketitle


\begin{abstract}
The Ekman boundary layer, a fundamental concept in geophysical fluid mechanics, describes the near-boundary fluid motion subject to rotation.
Within the singular limit framework of rapid rotation and vanishing viscosity, classical studies of Ekman theory (e.g., Desjardins and Grenier (1999), Masmoudi (2000)) are predominantly restricted to flat or small-amplitude boundary assumptions.

The conventional flat-boundary assumption obscures the complex mechanisms induced by topographic curvature; moreover, even small-amplitude perturbations reduce topographic effects to simple linear forcing terms.
Consequently, this paper investigates the singular limit behavior of rotating fluids over a non-flat boundary $z=B(x,y)$ of $\mathcal{O}(1)$ amplitude with uniformly bounded slope and curvature.
We elucidate how such topography modulates fluid dissipation through two distinct mechanisms: macroscopic topographic forcing and microscopic anisotropic pumping.
Specifically, at the macroscopic scale, large-scale topography forces the flow field to undergo kinematic uplift or subsidence.
At the microscopic scale, the local normal vector inclination breaks down the constant-thickness assumption of the classical Ekman boundary layer, leading to a spatially inhomogeneous boundary layer structure and inducing a novel generalized anisotropic Ekman pumping effect.

First, using multi-scale asymptotic analysis, we construct a class of approximate solutions that explicitly depend on the boundary's geometric characteristics, yielding a two-dimensional limit system fundamentally distinct from classical models.
A key innovation of this system is the introduction of a generalized velocity field defined via the topographic metric tensor.
This formulation not only generalizes the traditional isotropic linear damping to anisotropic geometric damping but also couples rotational effects to macroscopic vertical acceleration.
Furthermore, using energy methods, we establish the $L^2$ convergence of these variable-thickness approximate solutions to the weak solutions of the original three-dimensional system.
Finally, we analyze the multiple mechanisms governing rotating fluid motion over large-amplitude topography using a representative class of boundary geometries.
\end{abstract}



\section{Introduction}\label{D-sec1}
Rotating fluids constitute one of the core subjects in geophysical fluid dynamics and find extensive applications in the modeling of ocean circulation and atmospheric cyclones.
From a mathematical perspective, rotating fluids are typically governed by the following Navier-Stokes-Coriolis system:
\begin{equation}\label{D-1.1}
\begin{cases}
\partial_t \boldsymbol{u}^\varepsilon -\nu_{h}\Delta_{h} \boldsymbol{u}^\varepsilon-\nu_{v}\partial_z^2 \boldsymbol{u}^\varepsilon + (\boldsymbol{u}^\varepsilon \cdot \nabla)\boldsymbol{u}^\varepsilon + \varepsilon^{-1}\boldsymbol{R}\boldsymbol{u}^\varepsilon + \varepsilon^{-1}\nabla p^\varepsilon = 0, \\
\nabla \cdot \boldsymbol{u}^\varepsilon = 0,\\
\boldsymbol{u}^\varepsilon|_{t=0}=\boldsymbol{u}^\varepsilon_0,\quad
\boldsymbol{u}^\varepsilon|_{\partial\Omega}=0,
\end{cases}
\end{equation}
where $(t, \boldsymbol{x}) \in \mathbb{R}_+ \times \Omega$ with $\Omega\subset \mathbb{R}^3$. Here, $\boldsymbol{u}^\varepsilon$ and $p^\varepsilon$ denote the velocity and pressure of the fluid, respectively, whereas $\nu_h=\nu_h(\varepsilon)>0$ and $\nu_v=\nu_v(\varepsilon)>0$ describe the horizontal and vertical viscosity coefficients, respectively. The term $\varepsilon^{-1}\boldsymbol{R}\boldsymbol{u}^\varepsilon$ represents the Coriolis force, which accounts for the rotational effects of the system and takes the following form:
\begin{displaymath}
    \boldsymbol{R}\boldsymbol{u}^\varepsilon=\begin{pmatrix}
        0&-1&0\\
        1&0&0\\
        0&0&0
    \end{pmatrix}\boldsymbol{u}^\varepsilon. 
\end{displaymath}

When the fluid domain is the whole space $\mathbb{R}^3$ or the periodic domain $\mathbb{T}^3$, the well-posedness and vanishing viscosity limit of system \eqref{D-1.1} have been extensively studied. For instance, in $\Omega = \mathbb{T}^3$, Babin, Mahalov, and Nicolaenko\cite{Babin1996,Babin1999} demonstrated the global regularity of solutions to system \eqref{D-1.1} for sufficiently small $\varepsilon > 0$. For $\Omega = \mathbb{R}^3$, Chemin et al.\cite{Chemin2006} established unique global solution converging to a two-dimensional (2D) system as $\varepsilon \to 0$. 
These results indicate that, in the absence of boundaries, strong rotational effects significantly suppress three-dimensional (3D) vortex stretching, thereby degenerating the macroscopic flow field into a 2D columnar motion (the Taylor-Proudman effect). For further related results, we refer to Refs.~\refcite{Gallagher1998,Giga2005,Giga2008} and the references therein.

However, when the domain $\Omega$ has solid boundaries with Dirichlet no-slip conditions, the asymptotic behavior of the system differs fundamentally. Since the interior flow cannot directly accommodate the no-slip constraint at the boundary, the velocity gradient increases sharply within an extremely thin region near the boundary, thereby forming a distinct singular perturbation region—the Ekman boundary layer.

As a core concept in geophysical fluid dynamics, the Ekman boundary layer was first proposed in 1905 by the Swedish oceanographer V. W. Ekman\cite{Ekman1905}. This region is characterized by the balance between viscous friction, the Coriolis force, and the pressure gradient force, fundamentally governing the transport of momentum, heat, and mass between different fluid layers.

In rotating fluid models, the interaction between the Ekman boundary layer and solid boundaries gives rise to the classical Ekman pumping effect (see Refs. \refcite{Greenspan1968,Pedlosky}). Specifically, the no-slip condition disrupts the quasi-geostrophic balance; the induced secondary flow, subject to the incompressibility constraint, generates a vertical velocity at the outer edge of the boundary layer. Acting as either a source or a dissipation term within the macroscopic limit equations, this vertical flux couples the local boundary layer effects with the large-scale circulation.

At the theoretical analysis level, researchers primarily focus on characterizing the coupling mechanisms between the microscopic boundary layer and macroscopic dynamics. To handle such singular perturbation problems, the fluid velocity field $\boldsymbol{u}^\varepsilon$ is typically represented by the following multi-scale asymptotic expansion:
\begin{equation}\label{D-1.2}
\boldsymbol{u}^\varepsilon(t, \boldsymbol{x}) \sim \sum_{i=0}^{\infty} \delta^i \Big[ \boldsymbol{u}^{i,\text{I}}(t, x, y, z) + \boldsymbol{u}^{i,\text{BL}}\Big(t, x, y, \frac{d(\boldsymbol{x}, \partial\Omega)}{\delta}\Big) \Big],
\end{equation}
where $\boldsymbol{u}^{i,\text{I}}$ and $\boldsymbol{u}^{i,\text{BL}}$ denote the interior macroscopic flow field and the boundary layer corrector, respectively; $d(\boldsymbol{x}, \partial\Omega)$ represents the distance from a spatial point to the boundary, and $\delta$ is the characteristic thickness of the boundary layer. Within this theoretical framework, three fundamental issues need to be addressed: determining the boundary layer thickness $\delta$, solving for the structure of the boundary layer profile $\boldsymbol{u}^{i,\text{BL}}$, and establishing the convergence of this asymptotic expansion in appropriate function spaces.

However, classical Ekman boundary layer theories are predominantly confined to flat boundaries. When actual fluid domains feature complex, non-flat topographic characteristics (such as coastal regions, oceanic ridges, or trenches), the dynamic behavior of the Ekman boundary layer differs fundamentally from that over flat regions. This discrepancy arises primarily because the geometric features of the boundary significantly alter the distribution of the flow field and introduce additional complexities. Therefore, quantifying the impact of boundary geometric features on the Ekman boundary layer, the pumping effect, and the macroscopic flow structure is of profound significance for refining non-flat boundary layer theories and predicting flows in the vicinity of complex boundaries.

\subsection{Ekman Layer Theory over Flat Boundaries}\label{D-subsec1.1}

For a flat parallel-plate domain (i.e., $\Omega=\mathbb{T}^2 \times [0, 1]$ or $\Omega=\mathbb{R}^2 \times [0, 1]$), the boundary layer thickness of system \eqref{D-1.1} is the constant $\varepsilon$. Within this flat domain, addressing both well-prepared data and ill-prepared data cases, Grenier and Masmoudi\cite{Grenier1997,Masmoudi1998} constructed multi-scale approximate solutions for the Ekman boundary layer and established their convergence in the $L^2$ space.

To contrast the structural differences between flat and non-flat boundaries, we briefly recall the results of Grenier and Masmoudi\cite{Grenier1997} for the flat case: in the singular limit of vanishing viscosity, to overcome the nonlinear instability induced by the lack of dissipation and to ensure that the 3D weak solutions converge to the 2D limit system in the $L^2$ sense, the limit velocity field must satisfy a specific $L^\infty$ constraint. The specific theorem is as follows:
\begin{theorem}[Grenier-Masmoudi\cite{Grenier1997}]\label{GM1997}
Let $\Omega = \mathbb{T}^2 \times [0,1]$. Given an initial velocity $\bar{\boldsymbol{u}}_{0,h} \in H^4(\mathbb{T}^2)$, let $\bar{\boldsymbol{u}} = (\bar{\boldsymbol{u}}_h, 0)^T$ be the global strong solution to the following 2D Euler equations with damping:
\begin{equation*}
        \begin{cases}
            \partial_{t} \bar{\boldsymbol{u}}_h + (\bar{\boldsymbol{u}}_h \cdot \nabla_h) \bar{\boldsymbol{u}}_h + \sqrt{2} \bar{\boldsymbol{u}}_h + \nabla_h\bar{p} = 0, \\
            \nabla_h \cdot \bar{\boldsymbol{u}}_h = 0, \\
            \bar{\boldsymbol{u}}_h|_{t=0} = \bar{\boldsymbol{u}}_{0,h},
        \end{cases}
\end{equation*}
where $\nabla_h = (\partial_x, \partial_y)^T$ is the horizontal gradient operator.
Assume the initial data $\{\boldsymbol{u}_0^\varepsilon\}_{\varepsilon>0}\in L^2(\Omega)$, and let $\boldsymbol{u}^\varepsilon$ be the corresponding weak solution to system \eqref{D-1.1} with isotropic viscosity $\nu_h \sim \nu_v \sim \mathcal{O}(\varepsilon)$.

If $\lim_{\varepsilon \to 0} \lVert \boldsymbol{u}_0^\varepsilon - \bar{\boldsymbol{u}}_0 \rVert_{L^2(\Omega)} = 0$, and there exists a constant $C_0 > 0$ such that
\begin{equation}\label{D-1.3}
\lVert \bar{\boldsymbol{u}}_h \rVert_{L^\infty(\mathbb{T}^2)} \leqslant C_0,
\end{equation}
then as $\varepsilon \to 0$, we have
\begin{equation*}
\sup_{t \geqslant 0}\lVert \boldsymbol{u}^\varepsilon(t) - \bar{\boldsymbol{u}}(t) \rVert_{L^2(\Omega)}\rightarrow 0.
\end{equation*}
\end{theorem}
\begin{figure}[htbp]
	\centering
	\includegraphics[width=.8\textwidth]{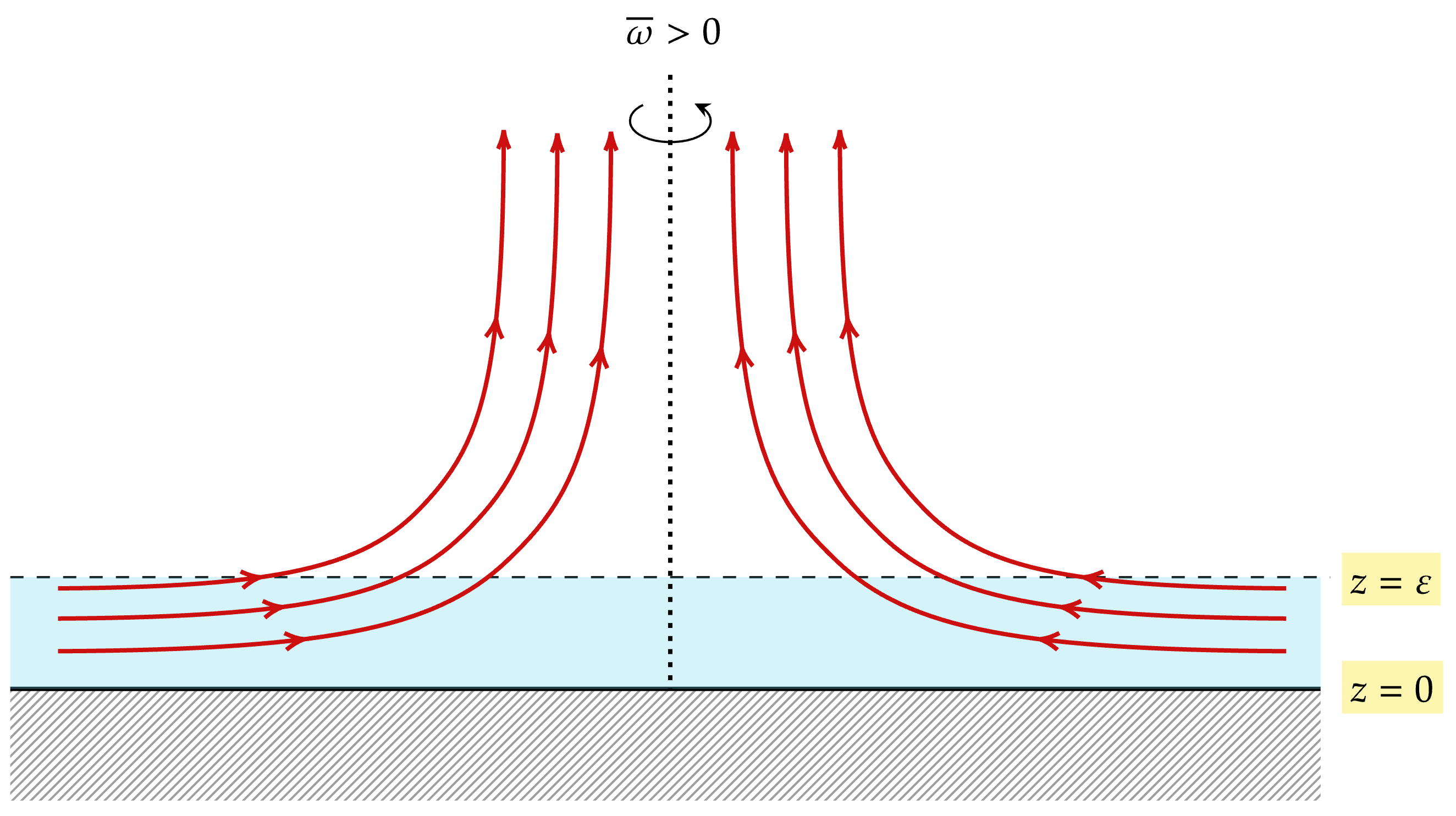}
	\caption{Ekman pumping induced by macroscopic positive vorticity over a flat boundary}\label{D-fig3}
\end{figure}

It should be emphasized that the domains $\mathbb{T}^2 \times[0, 1]$ and $\mathbb{R}^2 \times [0, 1]$ are a special class of geometric domains: their boundary unit normal vectors are globally constant ($\boldsymbol{n} = (0,0,1)^T$), and the vertical velocity of the limit flow field vanishes ($\bar{u}_3 = 0$). Under this flat-boundary assumption, Ekman pumping is dictated by the macroscopic vorticity $\bar{\omega}=\nabla_h^\bot\cdot\bar{\boldsymbol{u}}_h$ and degenerates into an isotropic linear damping term $\sqrt{2} \bar{\boldsymbol{u}}_h$ within the limit system (see Fig.\ref{D-fig3}). For further research developments concerning boundary layer theories in flat domains, we refer the reader to Refs. \refcite{Chemin2006,Gong2015,Liu2017,Liu2019,Liu20191,Liu2023}. Nevertheless, a systematic theoretical framework for the more general curved Ekman boundary layer problems has yet to be established.

\subsection{Ekman Layer Theory over Non-flat Boundaries}\label{D-subsec1.2}

Unlike the flat case, the thickness of a non-flat Ekman boundary layer depends on the boundary geometry, which poses significant analytical challenges. Even for certain specific curved geometries, the boundary layer often exhibits highly complex spatial structures.
For example, Stewartson\cite{Stewartson1957,Stewartson1966} discovered that in concentric rotating spheres, multiple nested sub-boundary layers with scales of $\mathcal{O}(\varepsilon^{1/2})$, $\mathcal{O}(\varepsilon^{1/3})$, and $\mathcal{O}(\varepsilon^{1/4})$ are generated near the inner boundary. Within special geometric domains such as cylindrical and spherical shells, Bresch et al.\cite{Bresch2004} and Rousset\cite{Rousset2007} constructed the corresponding approximate Ekman layer structures, respectively. The inherent difficulties lie in constructing approximate solutions at the intersection of the cylindrical sidewall and the horizontal boundary, as well as addressing the geometric constraints imposed by the equatorial exterior boundary on the boundary layer profile of spherical domains.

Beyond the aforementioned special geometries, for general non-flat boundaries, the complex nonlinear coupling between topographic undulations and the Ekman pumping mechanism makes the construction of boundary layers exceptionally challenging. To overcome this difficulty, previous studies have typically introduced small-amplitude or periodic boundary assumptions. For instance, when the horizontal boundary is periodic and of small amplitude, G\'erard-Varet and Bresch\cite{Bresch2005,Gerard2003} constructed approximate solutions for the Ekman boundary layer by introducing artificial interfaces.
Furthermore, in deriving wind-driven ocean circulation models, Desjardins and Grenier\cite{Desjardins1999} considered topography as an $\mathcal{O}(\varepsilon)$ perturbation (i.e., $z \in[-1+\varepsilon  B(x,y), 0]$) and derived the following 2D limit equations:
\begin{equation}\label{D-1.4}
    \begin{cases}
        \partial_t \bar{\boldsymbol{u}}_h + (\bar{\boldsymbol{u}}_h \cdot \nabla_h) \bar{\boldsymbol{u}}_h -  \Delta_h \bar{\boldsymbol{u}}_h + \frac{r_0}{2} \bar{\boldsymbol{u}}_h - (B + \beta y) \bar{\boldsymbol{u}}_h^\perp - \beta \boldsymbol{\tau} + \nabla_h \bar{p} = 0, \\
        \nabla_h \cdot \bar{\boldsymbol{u}}_h = 0,
    \end{cases} 
\end{equation}
where $\bar{\boldsymbol{u}}_h^\perp = (-\bar{u}_2, \bar{u}_1)^T$, $\boldsymbol{\tau}$ is the sea surface wind stress field, and $r_0$ and $\beta$ represent the Ekman pumping coefficient and the $\beta$-plane approximation parameter, respectively. 
As illustrated in system \eqref{D-1.4}, $\frac{r_0}{2} \bar{\boldsymbol{u}}_h$ denotes the linear damping term dominated by Ekman pumping, while $-B \bar{\boldsymbol{u}}_h^\perp$ represents the geometric effect induced by the small-amplitude topography. 
Similarly, Masmoudi\cite{Masmoudi2000} explored corresponding small-amplitude boundary layer problems within the domain $\mathbb{T}^2 \times[\varepsilon B(x,y), 1]$.
These small-amplitude models laid the foundation for understanding topographic effects on rotating flows, but they cannot capture the full geometric coupling observed in large-undulation terrains, which is closely related to the dynamics of topographic Rossby waves\cite{Rhines1973}.

Recently, Chemin et al.\cite{Chemin2024} investigated a special class of non-flat domains ($z \in [-B(x,y), 0]$). The topography considered in their study satisfies a specific geometric constraint where the isobaths are parallel to the coastline (i.e., $|\nabla_h B|^2 = F(B)$). However, this topographic assumption eliminates the component of the principal flow normal to the isobaths, which causes the nonlinear convective term to vanish during the limit process. Consequently, the limit system degenerates into an ordinary differential equation governed solely by the topographic Ekman pumping parameter $\lambda_B$, namely $\partial_t \bar{\bm{u}}_h+\lambda_B\bar{\bm{u}}_h=0$. For further relevant research concerning non-flat boundaries, we refer the reader to Refs. \refcite{Gerard2006,Gerard2010,Liu2014,Lopez2023,Ngo2020}.

The aforementioned non-flat boundary cases are all restricted by the geometric assumption of a nearly vertical boundary normal vector. Consequently, the thickness of the Ekman boundary layer degenerates to a constant $\varepsilon$, thereby obscuring the complex mechanisms induced by geometric characteristics over large-amplitude topography.

To address this limitation, this paper introduces an $\mathcal{O}(1)$ topography $z = B(x,y)$ and elucidates the effects of large-amplitude topography on rotating fluids from two distinct scales. First, the impermeability condition drives a macroscopic kinematic uplift or subsidence of the entire columnar flow along the topography. Second, the local normal inclination within the boundary layer modifies both the effective rotation rate and the boundary layer thickness—directly shifting the Ekman pumping capacity.
These mechanisms are closely related to physical processes near oceanic ridges, trenches, seamounts, and continental slopes. For example, the bottom slope can either enhance or suppress vertical transport. Local upwelling or downwelling may emerge near seamounts, while topographic curvature alters the interactions between vortices and the bottom. Therefore, it is essential to analyze the asymptotic behavior of these variable-thickness boundary layers and their coupling with geometric features like topographic inclination and principal curvature.

To better illustrate these topographic effects, we introduce two geometric characteristics of the boundary: the topographic inclination $\gamma$ (the angle between the local normal vector and the vertical) and the principal curvature.
\begin{itemize}
    \item The topographic inclination captures the first-derivative properties of the topography. The resulting geometric factor $\cos\gamma \in (0,1]$ reflects the steepness of the slope. Simply put, a smaller $\cos\gamma$ means steeper topography, making the boundary layer more susceptible to nonlinear separation.
    \item The principal curvature, a second-derivative property, describes the bending of the topographic surface. Its sign indicates whether the topography is concave or convex (like a basin versus a seamount), while its absolute value dictates the sharpness of the local topography.
\end{itemize}

Next, we consider the 3D Navier-Stokes-Coriolis system defined over the non-flat domain $\Omega=\mathbb{R}^2\times[B(x,y),B(x,y)+2]$:
\begin{equation}\label{D-1.5}
    \begin{cases}
        \partial_t \boldsymbol{u}^\varepsilon-\nu\varepsilon \Delta \boldsymbol{u}^\varepsilon+(\boldsymbol{u}^\varepsilon \cdot \nabla)\boldsymbol{u}^\varepsilon+\varepsilon^{-1}\boldsymbol{R}\boldsymbol{u}^\varepsilon+\varepsilon^{-1}\nabla p^\varepsilon=0,\\
        \nabla\cdot \boldsymbol{u}^\varepsilon=0,   \\
        \boldsymbol{u}^\varepsilon(t,x,y,z)|_{t=0}=\boldsymbol{u}^\varepsilon_0(x,y,z),\quad\boldsymbol{u}^\varepsilon(t,x,y,z)|_{\partial\Omega}=0,
    \end{cases}
\end{equation}
where $\nu\varepsilon>0$ denotes the viscosity coefficient, and the boundary is given by a smooth function $B(x,y)$. For simplicity in theoretical derivation, we assume the top and bottom boundaries are defined by parallel translations (i.e., the top boundary is $z=B(x,y)+2$).

For the dynamical system \eqref{D-1.5} with initial data $\boldsymbol{u}^\varepsilon_{0}\in {L}^2(\Omega)$, the skew-symmetric structure of the Coriolis operator $\varepsilon^{-1}\boldsymbol{R}\boldsymbol{u}^\varepsilon$ guarantees the existence of a global Leray-type weak solution $\boldsymbol{u}^\varepsilon$ under no-slip boundary conditions. By applying standard energy methods, this weak solution is shown to satisfy the following bound:
\begin{equation*}
    \lVert\boldsymbol{u}^\varepsilon\rVert_{{L}^\infty(\mathbb{R_+}; {L}^2(\Omega))}^2+\nu\varepsilon\int_0^\infty\lVert\nabla\boldsymbol{u}^\varepsilon(t,\cdot)\rVert_{{L}^2(\Omega)}^2\,dt   \leqslant \lVert\boldsymbol{u}^\varepsilon_{0}\rVert_{{L}^2(\Omega)}^2.
\end{equation*}

\subsection{Preliminaries and Mathematical Notation}\label{D-subsec1.3}

Before stating the main results of this paper, we introduce the key notation required for the subsequent analysis and derivations. First, we define the scenario where the boundary satisfies $B(x,y) \equiv \text{const}$ as the ``flat case''; all other scenarios are collectively referred to as the ``non-flat case''.

Second, for the non-flat surface $z = B(x,y)$, its unit normal vector is denoted by $\boldsymbol{n}$, and its direction cosines are expressed as:
\begin{equation*}
\boldsymbol{n} = \left( \tfrac{-B_x}{\sqrt{1+B_x^2+B_y^2}}, \tfrac{-B_y}{\sqrt{1+B_x^2+B_y^2}}, \tfrac{1}{\sqrt{1+B_x^2+B_y^2}} \right) =: (\cos\alpha, \cos\beta, \cos\gamma).
\end{equation*}
The second-order derivatives of $B(x,y)$ are characterized by its Hessian matrix $\boldsymbol{H}$:
\begin{equation*}
    \boldsymbol{H} = \begin{pmatrix} B_{xx} & B_{xy} \\ B_{xy} & B_{yy} \end{pmatrix}.
\end{equation*}
Furthermore, we have
\[
\nabla_h \gamma
= \cos^3\gamma \, \sin^{-1}\gamma \, \boldsymbol{H}\nabla_h B.
\]

Next, we define the induced metric matrix of the surface, $\boldsymbol{H}_0$, and its inverse as follows:
\begin{equation*}
    \boldsymbol{H}_0 = \begin{pmatrix} 1+B_x^2 & B_xB_y \\ B_xB_y & 1+B_y^2 \end{pmatrix},\quad
    \boldsymbol{H}^{-1}_0 = \cos^2\gamma\begin{pmatrix} 1+B_y^2 & -B_xB_y \\ -B_xB_y & 1+B_x^2 \end{pmatrix}.
\end{equation*}
It is straightforward to verify that the matrix $\boldsymbol{H}_0$ is strictly positive definite, and its determinant satisfies $\det \boldsymbol{H}_0 = \cos^{-2}\gamma$. To simplify the expressions in subsequent matrix operations, let $\boldsymbol{E}$ denote the second-order identity matrix, and we introduce the following skew-symmetric characteristic matrix:
\begin{equation*}
    \boldsymbol{E}_1 = \begin{pmatrix} 0 & 1 \\ -1 & 0 \end{pmatrix}.
\end{equation*}

Based on the aforementioned definitions, the Gaussian curvature $K_G$ and the mean curvature $K_A$ of the surface $B(x,y)$ can be expressed, respectively, as:
\begin{align*}
        K_G =& (\det \boldsymbol{H}_0)^{-2} \det \boldsymbol{H}, \\
        K_A =& \frac{1}{2}(\det \boldsymbol{H}_0)^{-\frac{3}{2}} \big[ B_{yy}(1+B_x^2) - 2B_{xy}B_xB_y + B_{xx}(1+B_y^2) \big],
\end{align*}
and the principal curvatures of the surface are given by
\[k_1, k_2 = K_A \pm \sqrt{K_A^2 - K_G}.
\]

Finally, we define the derivative operators and vectors as follows: the 3D gradient operator is denoted by $\nabla = (\partial_x, \partial_y, \partial_z)^T$, the horizontal gradient operator by $\nabla_h = (\partial_x, \partial_y)^T$, and the corresponding orthogonal 2D rotational gradient operator by $\nabla_h^\perp = (-\partial_y, \partial_x)^T$. Accordingly, the 3D and 2D Laplace operators are defined as $\Delta = \partial_x^2 + \partial_y^2 + \partial_z^2$ and $\Delta_h = \partial_x^2 + \partial_y^2$, respectively. The horizontal velocity vector is defined as $\bm{u}_h = (u_1, u_2)^T$, and its orthogonal counterpart is denoted by $\bm{u}_h^\perp = (-u_2, u_1)^T$.

\subsection{Main Results}\label{D-subsec1.4}
Compared to the classical scenarios of flat boundaries and small-amplitude topography, the $\mathcal{O}(1)$ non-flat boundary not only induces a macroscopic vertical velocity component but also generates anisotropic geometric damping through the spatially varying thickness of the boundary layer. Consequently, the structure of the limit system is fundamentally altered.

We first present the construction of approximate solutions for the Ekman boundary layer that incorporate geometric characteristics.
\begin{theorem}\label{D-th1}
    As $\varepsilon\to 0$, assume that the initial velocity field $\boldsymbol{u}_0^\varepsilon$ converges in $L^2(\Omega)$ to $\bar{\boldsymbol{u}}_0(x,y) = (\bar{\boldsymbol{u}}_{0,h}, \nabla_h B\cdot \bar{\boldsymbol{u}}_{0,h})^T$, with $\bar{\boldsymbol{u}}_{0,h}\in H^3(\mathbb{R}^2)$. Suppose that the boundary profile satisfies $B(x,y)\in W^{4,\infty}(\mathbb{R}^2)$, and there exist constants $C_1 \in (0,1)$ and $M_0 > 0$ such that
    \begin{equation}\label{D-1.6}
        0 < C_1 \leqslant \cos\gamma \leqslant 1, \qquad \sup_{(x,y)\in\mathbb{R}^2} |\nabla_h B| = M_0 > 0,
    \end{equation}
    where $\gamma$ denotes the angle between the boundary normal vector and the vertical direction. Then, system \eqref{D-1.5} possesses an approximate solution $(\boldsymbol{u}_{app}^\varepsilon, p_{app}^\varepsilon)$ of the following form:
\begin{equation*}
    \begin{cases}
        \boldsymbol{u}_{app}^\varepsilon = \bar{\bm{u}}(t,x,y) + \displaystyle\sum_{i=0}^1 \delta^i \boldsymbol{u}^i(t,x,y,z,\bar z,\tilde z) + \boldsymbol{u}^c, \\
        p_{app}^\varepsilon = \displaystyle\sum_{i=0}^2 \delta^i p^i(t,x,y,z,\bar z,\tilde z) + p^c,
    \end{cases}
\end{equation*}
where the boundary layer thickness is given by $\delta(x,y) = \sqrt{\nu}\varepsilon\cos^{-3/2}\gamma$, and $(\boldsymbol{u}^c,p^c)$ represents the divergence corrector.

In particular, the leading-order term of the asymptotic expansion, $\bar{\boldsymbol{u}}=(\bar{\boldsymbol{u}}_h,\bar u_3)^T$, is determined by the topographic metric tensor $\bm{H}_0$ and the generalized velocity $\bar{\bm{U}}_h$:
\begin{equation}\label{D-1.7}
\bar{\boldsymbol{u}}_h = \bm{H}_0^{-1} \bar{\bm{U}}_h, \qquad \bar{u}_3 = \cos^{2}\gamma \nabla_h B\cdot \bar{\bm{U}}_h,
\end{equation}
meanwhile, the pair $(\bar{\boldsymbol{u}}, \bar{p})$ satisfies the following 2D limit system:
\begin{equation}\label{D-1.8}
    \begin{cases}
        \partial_t \bar{\boldsymbol{U}}_h + ( \bar{\boldsymbol{u}}_h \cdot \nabla_h) \bar{\boldsymbol{U}}_h + \sqrt{\frac{\nu}{2}} \cos^{-\frac{1}{2}}\gamma \, \bar{\boldsymbol{U}}_h \\
        \qquad\qquad 
        +\cos^{2}\gamma
        \big(
            \sqrt{\frac{\nu}{2}}\cos^{-\frac32}\gamma\bm{E}
            -\bar{u}_3  \boldsymbol{H} \boldsymbol{E}_1
        \big)\boldsymbol{H}_0 \bar{\boldsymbol{U}}_h^\bot
        + \nabla_h \bar{p} = \boldsymbol{0}, \\[1.5ex]
        \nabla_h \cdot \bar{\boldsymbol{u}}_h = 0, \quad \bar{\boldsymbol{u}}_h|_{t=0} = \bar{\boldsymbol{u}}_{0,h}.
    \end{cases}
\end{equation}
\end{theorem}

\begin{remark}\label{D-rm1}
In what follows, a term-by-term comparison is conducted between the large-undulation topography limiting system \eqref{D-1.7}–\eqref{D-1.8} derived in this work and the small-amplitude topography model \eqref{D-1.4}.
\begin{enumerate}
\item {Nonlinear convection term}: The term $\partial_t \bar{\bm{U}}_h + (\bar{\bm{u}}_h \cdot \nabla_h) \bar{\bm{U}}_h$ in \eqref{D-1.8} depicts the transport of the generalized velocity field over the non-flat surface, corresponding to the convective term $\partial_t \bar{\bm{u}}_h + (\bar{\bm{u}}_h \cdot \nabla_h) \bar{\bm{u}}_h$ in \eqref{D-1.4}.

\item {Ekman damping effect}: The expression $\sqrt{\frac{\nu}{2}} \cos^{-\frac{1}{2}}\gamma \, \bar{\bm{U}}_h$ in \eqref{D-1.8} describes the anisotropic Ekman damping induced by the topographic tensor $\bm{H}_0$, whose intensity is determined by $\cos^{-\frac{1}{2}}\gamma$. This corresponds to the isotropic linear damping term $\frac{r_0}{2} \bar{\bm{u}}_h$ in \eqref{D-1.4}.

\item {Boundary geometric effect}: The term $\cos^{2}\gamma\big( \sqrt{\frac{\nu}{2}} \cos^{-\frac{3}{2}}\gamma \boldsymbol{E} - \bar{u}_3  \bm{H} \bm{E}_1 \big) \bm{H}_0 \bar{\bm{U}}_h^\bot$ in \eqref{D-1.8} reflects the nonlinear coupling among the Coriolis force, the macroscopic vertical velocity, and the geometric characteristics of the boundary over large-amplitude topography. In the small-amplitude model \eqref{D-1.4}, this effect manifests merely as the linear term $-B \bar{\bm{u}}_h^\bot$.
\end{enumerate}
In summary, introducing the generalized velocity $\bar{\bm{U}}_h$ not only clarifies the distinct physical mechanisms at play during the limit process but also demonstrates that the small-amplitude model \eqref{D-1.4} is simply a degenerate form of our limit system \eqref{D-1.8}.
\end{remark}

\begin{remark}\label{D-rm2}
If the limit system \eqref{D-1.8} is recast in terms of $\bar{\boldsymbol{u}}_h$, it can be equivalently expressed as:
\begin{equation}\label{D-1.9}
    \begin{cases}
        \partial_t \bar{\boldsymbol{u}}_h + (\bar{\boldsymbol{u}}_h \cdot \nabla_h) \bar{\boldsymbol{u}}_h + \sqrt{\frac{\nu}{2}} \cos^{-\frac32}\gamma (\cos\gamma \boldsymbol{H}_0  - \boldsymbol{E}_1 )\bar{\boldsymbol{u}}_h + (D_t \bar{u}_3)\nabla_h B  + \nabla_h {\bar{p}} = \bm{0},  \\
        \bar u_3=\nabla_h B\cdot \bar{\boldsymbol{u}}_h
        ,\quad D_t \bar{u}_3= (\partial_t  + (\bar{\boldsymbol{u}}_h \cdot \nabla_h)) \bar{u}_3,\\
        \nabla_h\cdot\bar{\boldsymbol{u}}_h = 0.
    \end{cases}
\end{equation}
The aforementioned limit equations reveal two distinct mechanisms of vertical motion induced by large-amplitude topography during rapid rotation:

The first mechanism involves the macroscopic topographic forcing flow of order $\mathcal{O}(1)$. Constrained by the impermeability condition at the boundary, the large-scale flow field generates a macroscopic vertical velocity $\bar{u}_3 = \bar{\boldsymbol{u}}_h \cdot \nabla_h B$. The resulting macroscopic vertical acceleration induces an additional pressure gradient, which subsequently translates into the topographic feedback term $(D_t \bar{u}_3)\nabla_h B$ within the horizontal momentum equations.

The second mechanism is the generalized Ekman pumping of order $\mathcal{O}(\delta)$. In the near-boundary region, the topographic inclination modifies the effective angular velocity of rotation, thereby causing the thickness of the Ekman boundary layer to vary spatially. Simultaneously, frictional effects disrupt the geostrophic balance, outputting an $\mathcal{O}(\delta)$ vertical pumping velocity into the interior flow field: $\sqrt{\frac{\nu}{2}} \cos^{-3/2}\gamma \big( \cos\gamma \boldsymbol{H}_0  - \boldsymbol{E}_1  \big)\bar{\boldsymbol{u}}_h.$ This expression explicitly reflects the threefold control exerted by the non-flat topography on the Ekman pumping process:

\begin{itemize}
    \item The variable-thickness property of the boundary layer influences the local pumping intensity (i.e., via the coefficient factor $\cos^{-3/2}\gamma$).
    \item The classical isotropic Ekman damping is transformed into an anisotropic geometric damping governed by the local metric tensor $\boldsymbol{H}_0$.
    \item The term $\sqrt{\frac{\nu}{2}} \cos^{-3/2}\gamma \boldsymbol{E}_1 \bar{\boldsymbol{u}}_h$ indicates that under the rapid rotation limit, the system still retains a partial residual Coriolis effect. In flat boundary scenarios (e.g., Ref. \refcite{Grenier1997}), this skew-symmetric term can be absorbed by the pressure gradient term. However, over non-flat topography, owing to the spatial dependence of the geometric metric coefficients, this term generally cannot be absorbed into the pressure.
\end{itemize}

It should be noted that the macroscopic vertical velocity within the limit flow field in this paper induces a $z$-dependent pressure gradient. Nevertheless, during the asymptotic expansion process, this term does not alter the 2D distribution of the limit flow field; instead, it is balanced by the Coriolis force generated by the first-order interior flow. For detailed derivations and proofs, we refer the reader to Subsubsection \ref{D-subsubsec3.2.3}.
\end{remark}

Having constructed the multi-scale approximate solutions incorporating geometric properties, we present the corresponding convergence results below.
\begin{theorem}\label{D-th2}
    Let $\boldsymbol{u}^\varepsilon \in L^\infty(\mathbb{R}_+; L^2(\Omega))$ be a family of weak solutions to the 3D Navier-Stokes-Coriolis system \eqref{D-1.5} with initial data $\boldsymbol{u}_0^\varepsilon \in L^2(\Omega)$. Assume that the following conditions hold:
    \begin{itemize}
        \item[1.]{Initial data conditions:}
        
        Let $\boldsymbol{u}^\varepsilon_0 = ({\boldsymbol{u}}^\varepsilon_{0,h}, {u}^\varepsilon_{0,3})$ and $\bar{\boldsymbol{u}}_{0,h} = \int_{[B, B+2]} \boldsymbol{u}_{0,h}^\varepsilon\,dz$. As $\varepsilon \to 0$, the initial data converge in the $L^2(\Omega)$ to:
\begin{equation}\label{D-1.10}
    \lim_{\varepsilon\to 0} \boldsymbol{u}^\varepsilon_0 = (\bar{\boldsymbol{u}}_{0,h}, \nabla_h B(x,y)\cdot \bar{\boldsymbol{u}}_{0,h})^T =: \bar{\boldsymbol{u}}_0,      
\end{equation}
and the limiting horizontal velocity satisfies
    \begin{equation}\label{D-1.11}
        \lVert \bar{\boldsymbol{u}}_{0,h} \rVert_{L^\infty(\mathbb{R}^2)} < \tfrac{\sqrt{\nu}}{4}
    \cos^{\frac52}\gamma,
    \end{equation}
    
        \item[2.]{Geometric constraints:}
    
    Suppose that the boundary profile satisfies $B(x,y)\in W^{4,\infty}(\mathbb{R}^2)$, and there exist constants $C_1 \in (0,1)$ and $M_0 > 0$ such that
\begin{equation}\label{D-1.12}
        0 < C_1 \leqslant \cos\gamma \leqslant 1, \qquad \sup_{(x,y)\in\mathbb{R}^2} |\nabla_h B| = M_0 > 0,
    \end{equation}
where $\gamma$ denotes the angle between the boundary normal vector and the vertical direction.
\end{itemize}
Then, as $\varepsilon \to 0$, we have:
\begin{equation}\label{D-1.13}
    \lim_{\varepsilon \to 0} \|\boldsymbol{u}^\varepsilon - \bar{\boldsymbol{u}}\|_{L^\infty(\mathbb{R}_+; L^2(\Omega))} = 0,
\end{equation}
where the limit flow field $\bar{\boldsymbol{u}}(t,x,y)$ satisfies the 2D system \eqref{D-1.7}-\eqref{D-1.8} incorporating topographic feedback and generalized Ekman damping.
\end{theorem}
\begin{remark}\label{D-rm3}
    This paper aims to reveal the geometric effects within non-flat boundary layers, and thus the well-prepared data condition \eqref{D-1.10} is adopted. It should be emphasized that if ill-prepared data are considered, Poincar\'e waves on a fast time scale will be generated in the rotating system (see Ref. \refcite{Chemin2006}). The analysis of such high-frequency oscillations is left for future work.
\end{remark}
\begin{remark}\label{D-rm4}
During the limit process, both horizontal and vertical dissipation degenerate. To prevent the Ekman boundary layer from undergoing separation, we must impose an upper bound constraint on the macroscopic velocity field: $\lVert \bar{\boldsymbol{u}}_{0,h} \rVert_{L^\infty(\mathbb{R}^2)} < \frac{\sqrt{\nu}}{4} \cos^{\frac52}\gamma$. In contrast to Ref.\refcite{Grenier1997}, this condition is additionally restricted by geometric effects: steeper topography (smaller $\cos\gamma$) corresponds to a lower upper limit of flow velocity.
\end{remark}
\begin{remark}\label{D-rm5}
The boundary constraint \eqref{D-1.12} is introduced for the following reasons. The lower bound $C_1 > 0$ limits the maximum topography slope and prevents a vertical boundary. It keeps the near-wall flow dominated by the Ekman mechanism (see Ref.\refcite{Chemin2006}) and suppresses flow separation over steep topography. The upper bound $\cos\gamma \leqslant 1$ permits local topography extrema such as hilltops and valleys. The condition $\sup |\nabla_h B| > 0$ excludes the globally flat boundary, so the system will not reduce to a flat or small-amplitude model. This ensures geometric effects dominate in the asymptotic expansion.
\end{remark}

The remainder of this paper is organized as follows. Section \ref{D-sec2} determines the thickness of the boundary layer by establishing a local boundary coordinate system and introducing a flattening transformation. Section \ref{D-sec3} uses multi-scale asymptotic analysis to obtain approximate solutions of system \eqref{D-1.5} in non-flat regions, and derives the limiting system with topographic feedback and Ekman pumping. Section \ref{D-sec4} proves the convergence between the approximate solutions and weak solutions of the system. Section \ref{D-sec5} analyzes multiple mechanisms governing rotating fluid flows over typical large-undulation topography. Finally, detailed proofs for propositions related to the 2D limiting system are presented in the appendix.

\section{Calculation of the Non-flat Ekman Boundary Layer Thickness}\label{D-sec2}

This section derives the thickness of the non-flat boundary layer. Unlike the classical flat model, our boundary exhibits $\mathcal{O}(1)$ spatial undulations, and hence its unit normal vector varies continuously with position. To address this, we use a local flattening transformation to derive the characteristic scale of the boundary layer.

First, denote the unit normal vector at an arbitrary point on the bottom boundary by $\boldsymbol{n} = (\cos\alpha, \cos\beta, \cos\gamma)$. We introduce local rotation angles $\theta$ and $\phi$ such that
\begin{equation}\label{D-2.1}
\begin{cases}
    \cos \theta = \pm \frac{\cos \gamma}{\sqrt{\cos^2 \alpha+\cos^2 \gamma}}, \\
    \sin \theta = \pm \frac{\cos \alpha}{\sqrt{\cos^2 \alpha+\cos^2 \gamma}},
\end{cases}
\quad
\begin{cases}
    \cos \phi = \pm \sqrt{\cos^2 \alpha + \cos^2 \gamma}, \\
    \sin \phi = \pm \cos \beta.
\end{cases}
\end{equation}

Using these definitions, we define the orthogonal rotation transformation between the local flattened coordinate system $\boldsymbol{x}' = (x', y', z')^T$ and the original Cartesian coordinate system $\boldsymbol{x} = (x, y, z)^T$:
\begin{equation}\label{D-2.2}
    \boldsymbol{x}' =
    \begin{pmatrix}
    x' \\ y' \\ z'
    \end{pmatrix}
    =
    \begin{pmatrix}
    \cos \theta & 0 & \sin \theta \\
    -\sin \theta \sin \phi & \cos \phi & \cos \theta \sin \phi \\
    -\sin \theta \cos \phi & -\sin \phi & \cos \theta \cos \phi
    \end{pmatrix}
    \begin{pmatrix}
    x \\ y \\ z
    \end{pmatrix}
    =: \boldsymbol{R}_0 \boldsymbol{x}. 
\end{equation}

\subsection{Linearized System in the Local Rotating Coordinate System}\label{D-subsec2.1}

\begin{figure}[htbp]\centering\includegraphics[width=.6\textwidth]{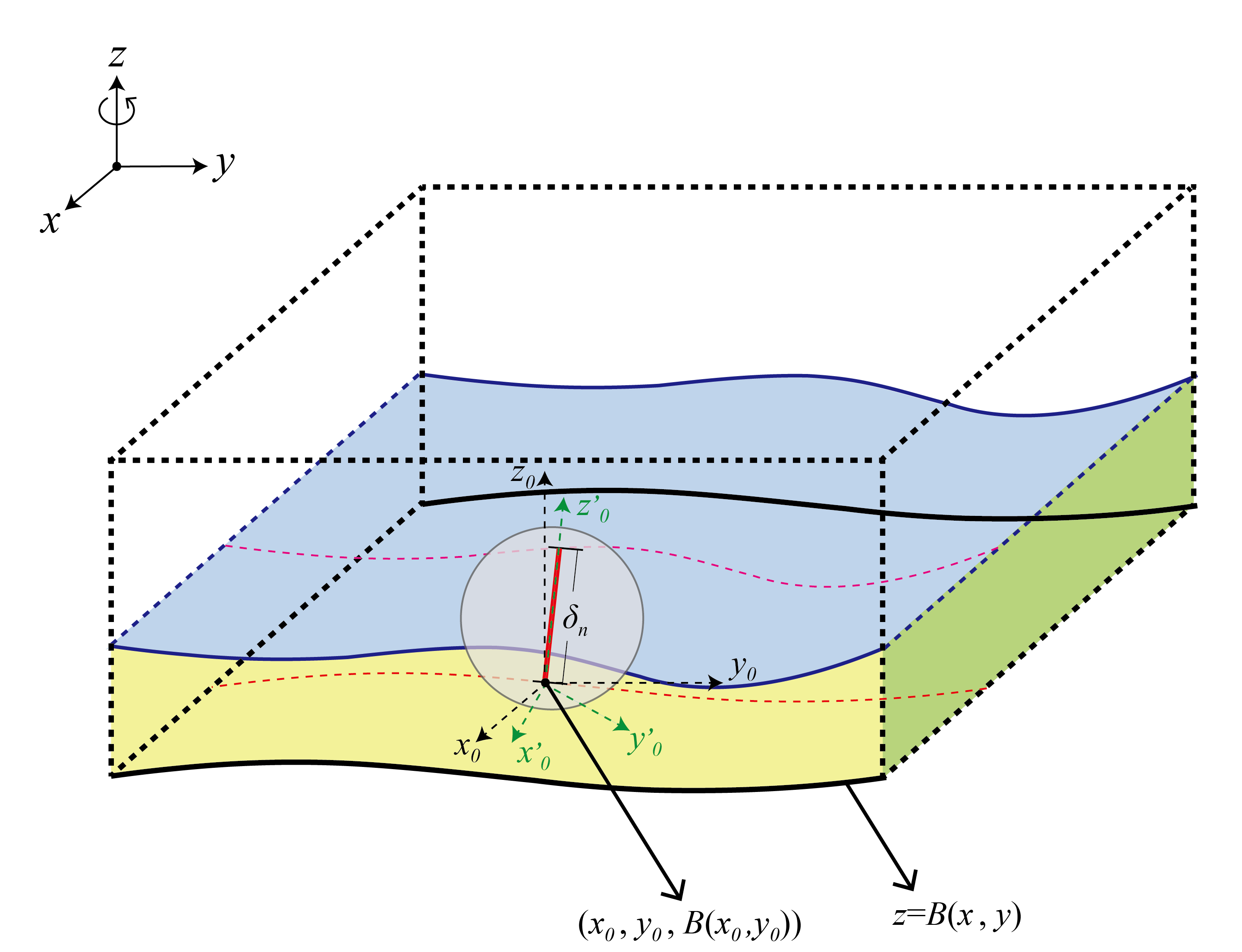}
    \caption{The boundary layer thickness at the local point $(x_0,y_0,B(x_0,y_0))$ on the bottom boundary.}\label{D-fig1}
\end{figure}

For an arbitrary point $(x_0, y_0)$ on the bottom boundary $z=B(x,y)$, the angle $\gamma$ between the normal vector $\boldsymbol{n}$ and the original $z$-axis varies spatially (see Fig.\ref{D-fig1}).

To derive the local characteristic scale of the boundary layer, we perform an analysis within a sufficiently small neighborhood of this point. In this local framework, the rotation parameters $\theta$ and $\phi$ can be taken as constants.

From the coordinate transformation \eqref{D-2.2}, the Jacobian matrix is given by:
\begin{equation}\label{D-2.3}
    \frac{\partial(x', y', z')}{\partial(x, y, z)} =
    \begin{pmatrix}
    \cos \theta & -\sin \theta \sin \phi & -\sin \theta \cos \phi \\
    0 & \cos \phi & -\sin \phi \\
    \sin \theta & \cos \theta \sin \phi & \cos \theta \cos \phi
    \end{pmatrix}
    =: \boldsymbol{R}_0^{-1}.
\end{equation}
The resulting spatial gradient operator is:
\begin{equation}\label{D-2.4}
    \nabla = 
    \begin{pmatrix}
    \partial_x \\ \partial_y \\ \partial_z
    \end{pmatrix}
    = \frac{\partial(x', y', z')}{\partial(x, y, z)}
    \begin{pmatrix}
    \partial_{x'} \\ \partial_{y'} \\ \partial_{z'}
    \end{pmatrix}
    =: \boldsymbol{R}_0^{-1} \nabla',
\end{equation}
Furthermore, by the properties of orthogonal transformations, the Laplace operator is form-invariant under this transformation, namely $\Delta = \partial_x^2 + \partial_y^2 + \partial_z^2 = \partial_{x'}^2 + \partial_{y'}^2 + \partial_{z'}^2 =: \Delta'$.

Note that in rotating fluid systems, the characteristic thickness of the boundary layer is fundamentally governed by the scale at which rotational effects balance viscous dissipation. Therefore, following the method from Ref.\refcite{Gerard2005}, we analyze the linearized system of \eqref{D-1.5}. Applying the local flattening transformation to this linearized system yields:
\begin{equation}\label{D-2.5}
\begin{cases}
    \partial_t \tilde{\boldsymbol{u}}^\varepsilon - \nu\varepsilon\Delta' \tilde{\boldsymbol{u}}^\varepsilon + \varepsilon^{-1} \tilde{\boldsymbol{R}} \tilde{\boldsymbol{u}}^\varepsilon + \varepsilon^{-1}\nabla' p^\varepsilon = 0, \\
    \nabla' \cdot \tilde{\boldsymbol{u}}^\varepsilon = 0,
\end{cases}
\end{equation}
where $\tilde{\boldsymbol{u}}^\varepsilon = \boldsymbol{R}_0 \boldsymbol{u}^\varepsilon$, and $\tilde{\boldsymbol{R}}$ is the equivalent Coriolis rotation matrix in the local rotating coordinate system:
\begin{equation}\label{D-2.6}
    \tilde{\boldsymbol{R}} :=
    \begin{pmatrix}
    0 & -\cos \theta \cos \phi & \cos \theta \sin \phi \\
    \cos \theta \cos \phi & 0 & -\sin \theta \\
    -\cos \theta \sin \phi & \sin \theta & 0
    \end{pmatrix}.
\end{equation}

\subsection{Calculation of the Characteristic Boundary Layer Thickness}\label{D-subsec2.2}

Performing a Laplace-Fourier analysis of the local linearized system \eqref{D-2.5}, we obtain
\begin{equation}\label{D-2.7}
    \begin{pmatrix}
    i\tau + \nu\varepsilon|\boldsymbol{\xi}|^2 & -\varepsilon^{-1}\cos\theta\cos\phi & \varepsilon^{-1}\cos\theta\sin\phi & \varepsilon^{-1}i\xi_1 \\
    \varepsilon^{-1}\cos\theta\cos\phi & i\tau + \nu\varepsilon|\boldsymbol{\xi}|^2 & -\varepsilon^{-1}\sin\theta & \varepsilon^{-1}i\xi_2 \\
    -\varepsilon^{-1}\cos\theta\sin\phi & \varepsilon^{-1}\sin\theta & i\tau + \nu\varepsilon|\boldsymbol{\xi}|^2 & \varepsilon^{-1}i\xi_3 \\
    i\xi_1 & i\xi_2 & i\xi_3 & 0
    \end{pmatrix}
    \begin{pmatrix}
    \hat{\tilde{u}}_1^\varepsilon \\ \hat{\tilde{u}}_2^\varepsilon \\ \hat{\tilde{u}}_3^\varepsilon \\ \hat{p}^\varepsilon
    \end{pmatrix}
    = 0.
\end{equation}
For this homogeneous system to possess non-trivial solutions, the determinant of the coefficient matrix must vanish. Since the geometric conditions \eqref{D-1.12} guarantee that the Ekman boundary layer dominates, it follows that the characteristic length $\delta_n$ of the bottom boundary layer along the local normal direction (i.e., the $z'$-axis direction) is:
\begin{equation}\label{D-2.8}
    \delta_{n} = \xi_3^{-1} = \varepsilon\sqrt{\nu} \cos^{-1/2}\gamma.
\end{equation}
\begin{figure}[htbp]
	\centering
	\includegraphics[width=.5\textwidth]{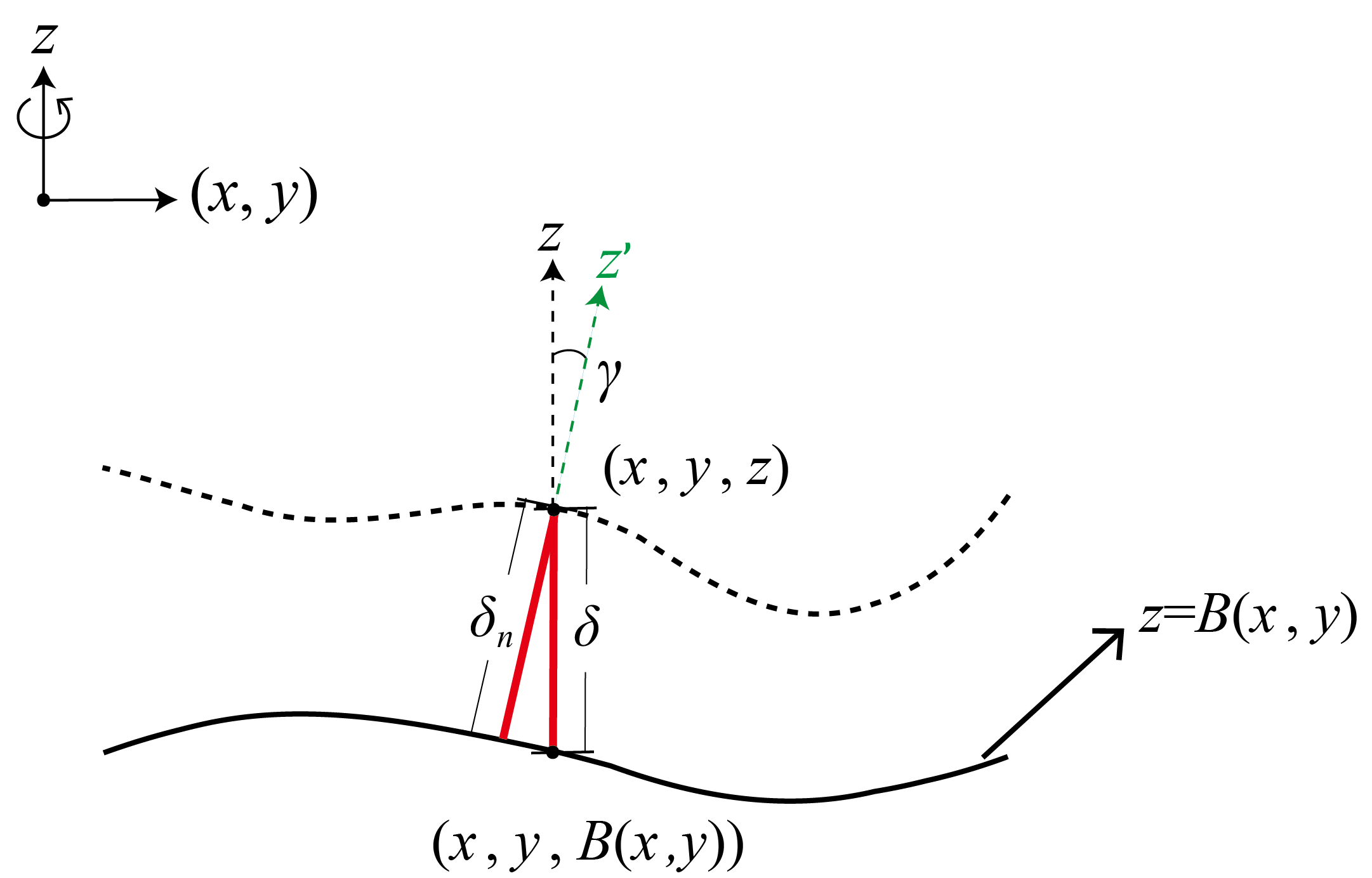}
	\caption{Effective boundary layer thickness $\delta$ in the original coordinate system.}\label{D-fig2}
\end{figure}

We then project the boundary layer thickness from the local normal direction onto the original $z$-axis direction (see Fig.\ref{D-fig2}). The effective boundary layer thickness $\delta(x,y)$ in the original coordinate system is thus given by:
\begin{equation}\label{D-2.9}
    \delta(x,y) = \frac{\delta_{n}(x,y)}{\cos \gamma(x,y)} = \varepsilon\sqrt{\nu} \cos^{-3/2}\gamma.
\end{equation}

\begin{remark}\label{D-rm6}
The above derivation indicates that steeper topographic implies a smaller $\cos\gamma$ and a thicker boundary layer. Specifically, large undulations of $\mathcal{O}(1)$ magnitude lead to spatial variation of the thickness $\delta$, which fundamentally distinguishes our asymptotic analysis from conventional models assuming constant boundary layer thickness. This reveals that large-undulation topographic not only imposes external boundary forcing, but also alters the intrinsic characteristic scale of the boundary layer.
\end{remark}

\section{Construction of Approximate Solutions}\label{D-sec3}

This section focuses on constructing multi-scale approximate solutions for the Navier-Stokes-Coriolis system \eqref{D-1.5} over non-flat domains. By asymptotically matching the interior flow with the boundary layer correctors, we apply a multi-scale expansion to the nonlinear system to obtain their explicit structure.

To ensure a smooth matching between the interior flow and the boundary layers, we introduce a cutoff function $\chi(z- B(x,y))\in C^\infty$ satisfying
\begin{equation}\label{D-3.1}
	\chi(z- B)=
	\begin{cases}
		\begin{aligned}
			&0,  && 0 \leqslant z-B \leqslant 1/2,\\
			&1,  && 3/2 \leqslant z-B \leqslant 2.
		\end{aligned}
	\end{cases}
\end{equation}
Evidently, the support of any $k$-th order derivative of $\chi$ ($k \ge 1$) lies strictly in the interval $(B+1/2, B+3/2)$.

We decompose the approximate solution $(\boldsymbol{u}_{\text{app}}^\varepsilon, p_{\text{app}}^\varepsilon)$ into the interior flow field (denoted by superscript $\mathrm{I}$) and the bottom and top boundary layer correctors (denoted by superscripts $\mathrm{b}$ and $\mathrm{t}$ respectively), namely
\begin{equation}\label{D-3.2}
\begin{cases}
    \boldsymbol{u}_{\text{app}}^\varepsilon  = \sum\limits_{i=0}^m\delta^i 
    \big(
    \boldsymbol{u}^{i,\mathrm{I}} + (1-\chi)\boldsymbol{u}^{i,\text{BL,b}}(t, \tilde x, \tilde y, \tilde{z}) +\chi \boldsymbol{u}^{i,\text{BL,t}}(t, \bar x, \bar y, \bar{z})   
    \big), \\
    p_{\text{app}}^\varepsilon  = \sum\limits_{i=0}^m \delta^i \big(p^{i,\mathrm{I}} + (1-\chi) p^{i,\text{BL,b}}(t, \tilde x, \tilde y, \tilde{z}) +\chi p^{i,\text{BL,t}}(t,\bar x, \bar y, \bar{z}) \big),
\end{cases}
\end{equation}
where the boundary layer thickness $\delta(x,y) = \sqrt{\nu\varepsilon} \cos^{-3/2}\gamma(x,y)$ varies spatially.

For both boundaries, we introduce the following fast variables:
\begin{equation*}
\tilde{z} = \frac{z-B(x,y)}{\delta(x,y)}, \qquad \bar{z} = \frac{B(x,y)+2-z}{\delta(x,y)}.
\end{equation*}
To ensure the validity of the multi-scale expansion, the interior flow and the boundary layers must satisfy the following matching conditions. 
On the one hand, to satisfy the Dirichlet boundary conditions, the boundary layer correctors must cancel the non-zero trace of the interior flow field, i.e.,
\begin{equation}\label{D-3.3}
\begin{cases}
    \boldsymbol{u}^{i,\text{BL,b}}(t,\tilde x,\tilde y, 0) = -\boldsymbol{u}^{i,\mathrm{I}}(t, x, y, z)\big|_{z=B(x,y)}, \\
    \boldsymbol{u}^{i,\text{BL,t}}(t,\bar x,\bar y, 0) = -\boldsymbol{u}^{i,\mathrm{I}}(t, x, y, z)\big|_{z=B(x,y)+2}.
\end{cases}
\end{equation}
On the other hand, these correctors must decay rapidly far from the boundaries. This gives the following asymptotic matching conditions:
\begin{equation}\label{D-3.4}
\lim_{\tilde{z}\to +\infty}\boldsymbol{u}^{i,\text{BL,b}}(t, \tilde x, \tilde y, \tilde{z}) = \boldsymbol{0}, \qquad \lim_{\bar{z}\to +\infty}\boldsymbol{u}^{i,\text{BL,t}}(t,\bar x, \bar y, \bar{z}) = \boldsymbol{0}.
\end{equation}

Note that by the support property of the cutoff function $\chi$, its spatial derivatives are non-zero only in the region $(B+1/2, B+3/2)$ away from the boundaries. Together with the boundary layer condition \eqref{D-3.4}, the coupling terms involving derivatives of $\chi$ may be neglected when constructing the local asymptotic expansion of the boundary layer.

As the two boundaries are parallel, the subsequent asymptotic derivation focuses on the bottom boundary layer.

\subsection{Boundary Layer Coordinate Transformations and Operator Expansions}\label{D-subsec3.1}

To characterize the flow structure within the boundary layers over non-flat case, we first introduce coordinate transformations for both boundary layers, namely:
\begin{equation*}
(\tilde{x}, \tilde{y}, \tilde{z}) = \big(x, y, \tfrac{z-B(x,y)}{\sqrt{\nu\varepsilon}}\cos^{\frac{3}{2}}\gamma\big), \quad (\bar{x}, \bar{y}, \bar{z}) = \big(x, y, \tfrac{B(x,y)+2-z}{\sqrt{\nu\varepsilon}}\cos^{\frac{3}{2}}\gamma\big).
\end{equation*}
When taking spatial derivatives under the above boundary layer coordinate system, the topographic undulations introduce geometric coupling terms and differential singularities.
To resolve this, we asymptotically decompose the gradient operator $\nabla$ and the Laplace operator $\Delta$, using the small parameter $\delta$ as the expansion scale.

For simplicity in what follows, we define the horizontal differential operators $\nabla_{\tilde{h}} = (\partial_{\tilde{x}}, \partial_{\tilde{y}})^T$ and $\Delta_{\tilde{h}} = \partial_{\tilde{x}}^2 + \partial_{\tilde{y}}^2$. For the bottom boundary layer, these operators can be expanded as follows:
\begin{align}
    \nabla &= \tilde{\nabla}_{\delta^0} + \delta^{-1}\tilde{\nabla}_{\delta^{-1}},\label{D-3.5}\\
    \Delta &= \tilde{\Delta}_{\delta^0} + \delta^{-1}\tilde{\Delta}_{\delta^{-1}} + \delta^{-2}\tilde{\Delta}_{\delta^{-2}},\label{D-3.6}
\end{align}
where the explicit expressions for the operators at each order are given by:
\begin{equation*}
\begin{aligned}
    \tilde{\nabla}_{\delta^0} =& \left(\nabla_{\tilde{h}} + \tfrac{3\tilde{z}}{2}\nabla_h(\ln\cos\gamma)\partial_{\tilde{z}}, \, 0\right)^T, \\[1mm]
    \tilde{\Delta}_{\delta^0} =& \Delta_{\tilde{h}} + \tfrac{9}{4}|\tilde{z}\,\nabla_h(\ln\cos\gamma)|^2\partial_{\tilde{z}}^2 + 3\tilde{z}\nabla_h(\ln\cos\gamma)\cdot\nabla_{\tilde{h}}\partial_{\tilde{z}} \\
    &+ \tfrac{3\tilde{z}}{2}\left(
        \tfrac{3}{2}|\nabla_h(\ln\cos\gamma)|^2
        +\Delta_h(\ln\cos\gamma)\right)\partial_{\tilde{z}}, \\[1mm]
    \tilde{\nabla}_{\delta^{-1}} =& \nabla(\delta\tilde{z})\partial_{\tilde{z}}, \\[1mm]
    \tilde{\Delta}_{\delta^{-1}} =& 2\nabla_h(\delta\tilde{z})\cdot\nabla_{\tilde{h}}\partial_{\tilde{z}} + 3\tilde{z}\nabla_h(\ln\cos\gamma)\cdot\nabla_h(\delta\tilde{z})\partial_{\tilde{z}}^2\\
    &+ \Delta_h(\delta\tilde{z})\partial_{\tilde{z}} + 3\nabla_h(\ln\cos\gamma)\cdot\nabla_h(\delta\tilde{z})\partial_{\tilde{z}}, \\[1mm]
    \tilde{\Delta}_{\delta^{-2}} =& \cos^{-2}\gamma \partial_{\tilde{z}}^2.
\end{aligned}
\end{equation*}

For the top boundary layer, the expansion structure of the differential operators is identical to that of the bottom boundary layer. To obtain the specific expressions, simply replace the local variables $(\tilde{x},\tilde{y},\tilde{z})$ with $(\bar{x},\bar{y},\bar{z})$. The decompositions are:
\begin{equation*}
    \nabla = \bar{\nabla}_{\delta^0} + \delta^{-1}\bar{\nabla}_{\delta^{-1}}, \qquad
    \Delta = \bar{\Delta}_{\delta^0} + \delta^{-1}\bar{\Delta}_{\delta^{-1}} + \delta^{-2}\bar{\Delta}_{\delta^{-2}}.
\end{equation*}

\subsection{Asymptotic Expansion Analysis}\label{D-subsec3.2}

Substituting the multi-scale expansion \eqref{D-3.2} into the original system \eqref{D-1.5}. Equating coefficients of like powers of $\delta^i$, we derive the governing equations order by order for $i=-2,-1,0,\dots$. To simplify notation, we let $\boldsymbol{u}^i$ denote the $i$-th order velocity profile, namely:
\begin{equation*}
    \boldsymbol{u}^i = \boldsymbol{u}^{i,\mathrm{I}} + (1-\chi)\boldsymbol{u}^{i,\mathrm{BL,b}} + \chi\boldsymbol{u}^{i,\mathrm{BL,t}}.
\end{equation*}
Similarly, the $i$-th order pressure term is given by $p^i = p^{i,\mathrm{I}} + p^{i,\mathrm{BL,b}} + p^{i,\mathrm{BL,t}}$.

\subsection{Asymptotic Expansion Analysis}\label{D-subsec3.2}

Substitute the multi-scale expansion \eqref{D-3.2} into the original system \eqref{D-1.5}. By equating coefficients of like powers of $\delta^i$, we derive the governing equations order by order for $i=-2,-1,0,\dots$. To simplify subsequent expressions, we denote the $i$-th order velocity profile uniformly by $\boldsymbol{u}^i$, namely:
\begin{equation*}
    \boldsymbol{u}^i = \boldsymbol{u}^{i,\mathrm{I}} + (1-\chi)\boldsymbol{u}^{i,\mathrm{BL,b}} + \chi\boldsymbol{u}^{i,\mathrm{BL,t}}.
\end{equation*}

\subsubsection{The $\mathcal{O}(\delta^{-2})$ Order Terms in the Asymptotic Expansion}\label{D-subsubsec3.2.1}

We first examine the $\mathcal{O}(\delta^{-2})$ order terms in the asymptotic expansion:
\begin{displaymath} 
	\partial_{\tilde{z}}p^{0,\mathrm{BL,b}}=\partial_{\bar{z}}p^{0,\mathrm{BL,t}}=0.
\end{displaymath}
Together with the far-field decay conditions $\lim_{\tilde{z}\to+\infty}p^{0,\mathrm{BL,b}}=\lim_{\bar{z}\to+\infty}p^{0,\mathrm{BL,t}}=0$, it follows that 
\[
p^{0,\mathrm{BL,b}}=p^{0,\mathrm{BL,t}}= 0.
\]
This implies that the zero-order boundary layer pressure correction terms are identically zero locally.

\subsubsection{The $\mathcal{O}(\delta^{-1})$ Order Terms in the Asymptotic Expansion}\label{D-subsubsec3.2.2}

We now analyze the $\mathcal{O}(\delta^{-1})$ order interior equations, which describe the geostrophic balance between the Coriolis force and pressure gradient:
\begin{equation}\label{D-3.7}
    \begin{cases}
        -\boldsymbol{E_1}\boldsymbol u_h^{0,\mathrm{I}}+\nabla_hp^{0,\mathrm{I}}=\boldsymbol{0},\\
        \partial_{z} p^{0,\mathrm{I}}=0.
    \end{cases}
\end{equation}
Combined with the incompressibility condition $\nabla \cdot\boldsymbol  u^{0,\mathrm{I}}=0$, it follows that both $\boldsymbol u^{0,\mathrm{I}}$ and $p^{0,\mathrm{I}}$ are independent of the vertical coordinate $z$. This manifests the Taylor-Proudman theorem: the leading-order interior flow has a 2D columnar structure. Applying the horizontal divergence operator to \eqref{D-3.7} yields the zero-order interior pressure:
\begin{equation}\label{D-3.8}
    p^{0,\mathrm{I}}=\Delta_h^{-1}\nabla^\perp_h \cdot\boldsymbol u_h^{0,\mathrm{I}}.
\end{equation}

We next analyze the $\mathcal{O}(\delta^{-1})$ order terms in the bottom boundary layer, where the contribution of nonlinear advection terms is $(u_3^0 - \boldsymbol{u}_h^0 \cdot \nabla_h B)\partial_{\tilde{z}}\boldsymbol{u}_h^{0,\mathrm{BL,b}} $. This gives:
\begin{equation}\label{D-3.9}
    \begin{cases}
        \cos\gamma\partial_{\tilde z}^2 \boldsymbol{u}_h^{0,\mathrm{BL,b}} + \boldsymbol{E}_1 \boldsymbol{u}_h^{0,\mathrm{BL,b}} + \nabla_h B \, \partial_{\tilde z}p^{1,\mathrm{BL,b}} + (u_3^0 - \boldsymbol{u}_h^0 \cdot \nabla_h B)\partial_{\tilde{z}}\boldsymbol{u}_h^{0,\mathrm{BL,b}} = \boldsymbol{0},\\[1ex]
        \partial_{\tilde z}p^{1,\mathrm{BL,b}} = \cos\gamma\partial_{\tilde z}^2 u_3^{0,\mathrm{BL,b}} + (u_3^0 - \boldsymbol{u}_h^0 \cdot \nabla_h B)\partial_{\tilde{z}}u_3^{0,\mathrm{BL,b}},\\[1ex]
        \partial_{\tilde z}u_3^{0,\mathrm{BL,b}} = \nabla_hB \cdot \partial_{\tilde z}\boldsymbol{u}_h^{0,\mathrm{BL,b}}.
    \end{cases}
\end{equation}

Integrating $\eqref{D-3.9}_3$ over $\tilde{z}$ with the far-field decay condition $\lim_{\tilde{z}\to+\infty}\boldsymbol{u}^{0,\mathrm{BL,b}}=\boldsymbol{0}$ gives the bottom boundary layer vertical velocity:
\begin{equation}\label{D-3.10}
    u_3^{0,\mathrm{BL,b}} = \nabla_h B \cdot \boldsymbol{u}_h^{0,\mathrm{BL,b}}.
\end{equation}
Similarly, the top boundary layer vertical velocity is
\begin{equation}\label{D-3.11}
    u_3^{0,\mathrm{BL,t}}  = \nabla_h B \cdot \boldsymbol{u}_h^{0,\mathrm{BL,t}}.
\end{equation}
Furthermore, combining the boundary matching conditions \eqref{D-3.3} with the $z$-independence of $\boldsymbol{u}^{0,\mathrm{I}}$, we obtain the leading-order interior vertical velocity:
\begin{equation}\label{D-3.12}
    u_3^{0,\mathrm{I}} = \nabla_h B \cdot \boldsymbol{u}_h^{0,\mathrm{I}}.
\end{equation}
\begin{remark}
    Equation \eqref{D-3.12} shows that subject to $\boldsymbol{u}^{0,\mathrm{I}} \cdot \boldsymbol{n} = 0$, the interior flow induces a vertical velocity $u_3^{0,\mathrm{I}}$. This term vanishes for flat or small-amplitude topography, but dominates vertical motion in our large-undulation setting: fluid climbs or descends along topography as it moves horizontally.
\end{remark}

Substituting \eqref{D-3.10}–\eqref{D-3.12} into \eqref{D-3.9} eliminates the nonlinear terms.
Eliminating the pressure gradient term in $\eqref{D-3.9}_1$ by $\eqref{D-3.9}_2$ and $\eqref{D-3.9}_3$ simplifies it to
\begin{equation*}   
    \boldsymbol{H_0} \partial_{\tilde z}^2\boldsymbol u_h^{0,\mathrm{BL,b}}+\cos^{-1}\gamma\boldsymbol{E_1}\boldsymbol u_h^{0,\mathrm{BL,b}}=\boldsymbol{0},
\end{equation*}
Combined with \eqref{D-3.3} and \eqref{D-3.4}, the zero-order bottom boundary layer velocity is:
\begin{equation}\label{D-3.13}
\boldsymbol{u}_h^{0,\mathrm{BL,b}}=
\boldsymbol{M}(\tfrac{\tilde z}{\sqrt 2}) \boldsymbol{u}_h^{0,\mathrm{I}},
\end{equation}
where the coefficient matrix is defined by
\begin{equation*}
	\boldsymbol{M}(Z)=
	-\mathrm{e}^{-Z}
	\Big(
	\cos(Z)\boldsymbol{E}+\sin(Z)\cos\gamma\boldsymbol{E_1}\boldsymbol{H_0}
	\Big)\,.
\end{equation*}
Equation \eqref{D-3.13} shows non-flat topography modifies both boundary layer thickness and Ekman spiral structure. In the flat case where $\boldsymbol{H}_0=\boldsymbol{E}$ and $\cos\gamma=1$, the classical rotational exponentially decaying structure is recovered. In the non-flat case, however, the term $\boldsymbol{E}_1\boldsymbol{H}_0$ subjects the rotation to local metric distortions, acting as the fundamental mechanism underlying the anisotropic boundary layer.

Integrating $\eqref{D-3.9}_2$ vertically gives the first-order bottom boundary layer correction pressure:
\begin{equation}\label{D-3.14}
p^{1,\mathrm{BL,b}} = \mathrm{e}^{-\frac{\tilde z}{\sqrt 2}} \Big( \sin(\tfrac{\tilde z}{\sqrt 2}+\tfrac\pi4)\cos\gamma\nabla_h B + \sin(\tfrac{\tilde z}{\sqrt 2}-\tfrac\pi4)\cos^2\gamma\nabla^\perp_h B\Big) \cdot \boldsymbol u_h^{0,\mathrm{I}}.
\end{equation}

Similarly, the top boundary layer velocity and pressure are:
\begin{align}
\boldsymbol{u}^{0,\mathrm{BL,t}}_h&=
\boldsymbol{M}(\tfrac{\bar z}{\sqrt 2}) \boldsymbol{u}_h^{0,\mathrm{I}}, \label{D-3.15} \\[3pt]
p^{1,\mathrm{BL,t}} 
&= -\mathrm{e}^{-\frac{\bar z}{\sqrt 2}} \Big( \sin(\tfrac{\bar z}{\sqrt 2}+\tfrac\pi4)\cos\gamma\nabla_h B + \sin(\tfrac{\bar z}{\sqrt 2}-\tfrac\pi4)\cos^2\gamma\nabla^\perp_h B\Big) \cdot \boldsymbol u_h^{0,\mathrm{I}}.\label{D-3.16}
\end{align}

\subsubsection{The $\mathcal{O}(1)$ Order Interior Terms in the Asymptotic Expansion}\label{D-subsubsec3.2.3}

The $\mathcal{O}(1)$ order interior equations of the asymptotic expansion are given by:
\begin{equation}\label{D-3.17}
    \begin{cases}
        \partial_t \boldsymbol{u}^{0,\mathrm{I}}_h + (\boldsymbol{u}^{0,\mathrm{I}}_h \cdot \nabla_h) \boldsymbol{u}^{0,\mathrm{I}}_h - \varepsilon^{-1}\delta \boldsymbol{E}_1 \boldsymbol{u}^{1,\mathrm{I}}_h + \nabla_h (\varepsilon^{-1}\delta {p}^{1,\mathrm{I}}) = \boldsymbol{0}, \\[1ex]
        \partial_t {u}^{0,\mathrm{I}}_3 + (\boldsymbol{u}^{0,\mathrm{I}}_h \cdot \nabla_h) {u}^{0,\mathrm{I}}_3 + \partial_z (\varepsilon^{-1}\delta {p}^{1,\mathrm{I}}) = 0, \\[1ex]
        \nabla_h \cdot \boldsymbol{u}^{0,\mathrm{I}}_h = 0.
    \end{cases}
\end{equation}
First, integrating the vertical momentum equation $\eqref{D-3.17}_2$ along the $z$-direction yields the first-order interior pressure field:
\begin{equation}\label{D-3.22}
    \varepsilon^{-1}\delta {p}^{1,\mathrm{I}}(t,x,y,z) = -z D_t {u}^{0,\mathrm{I}}_3 + {\bar{p}}(t,x,y),
\end{equation}
where $D_t = \partial_t + {\boldsymbol{u}}^{0,\mathrm{I}}_h \cdot \nabla_h$. Substituting this pressure expression back into the horizontal momentum equation and differentiating with respect to $z$ gives the following relation:
\begin{equation}\label{D-3.23}
    \varepsilon^{-1}\delta \partial_z \boldsymbol{u}_h^{1,\mathrm{I}} =  \nabla_h^\perp \big(\partial_z (\varepsilon^{-1}\delta {p}^{1,\mathrm{I}})\big) = -\nabla_h^\perp \big( D_t {u}^{0,\mathrm{I}}_3 \big).
\end{equation}
This shows that the geometry-induced pressure gradient is balanced by the vertical shear of the first-order interior horizontal velocity.

Integrating equation \eqref{D-3.23} along the vertical direction gives the structure of the first-order horizontal velocity:
\begin{equation}\label{D-3.24}
    \varepsilon^{-1}\delta \boldsymbol{u}_h^{1,\mathrm{I}}(t,x,y,z) = - z \nabla_h^\perp \big(D_t {u}^{0,\mathrm{I}}_3 \big) + {\bar{\boldsymbol{u}}}_h^{1}(t,x,y),
\end{equation}
where the structure of the 2D component ${\bar{\boldsymbol{u}}}_h^{1}$ is to be determined.

Substituting equations \eqref{D-3.22} and \eqref{D-3.24} back into the original horizontal momentum equation, equation $\eqref{D-3.17}_1$ can be rewritten as:
\begin{equation}\label{D-3.25}
    \partial_t {\boldsymbol{u}}_h^{0,\mathrm{I}} + ({\boldsymbol{u}}^{0,\mathrm{I}}_h \cdot \nabla_h) {\boldsymbol{u}}^{0,\mathrm{I}}_h - \boldsymbol{E}_1 {\bar{\boldsymbol{u}}}_h^1(t,x,y) + \nabla_h {\bar{p}}(t,x,y) = \boldsymbol{0}.
\end{equation}

Furthermore, to determine the form of ${\bar{\boldsymbol{u}}}_h^1$ and close system \eqref{D-3.25}, we use the incompressibility condition for the $\mathcal{O}(1)$ order boundary layers, i.e.,
\begin{align}
	\bm{u}^{1,\mathrm{BL,b}}\cdot\bm{n}
	=&\cos\gamma	
	\int_{\tilde{z}}^{+\infty}
	\tilde{\nabla}_{\delta^0}
	\cdot\boldsymbol{u}^{0,\mathrm{BL,b}}
	d\tilde{z}'\,,\label{D-3.18}\\
    \bm{u}^{1,\mathrm{BL,t}}\cdot\bm{n}
	=&
    -\cos\gamma\int_{\bar{z}}^{+\infty}
	\bar{\nabla}_{\delta^0}
	\cdot\boldsymbol{u}^{0,\mathrm{BL,t}}
	d\bar{z}'\,.\label{D-3.19}
\end{align}
Combined with the boundary matching conditions \eqref{D-3.3}, we derive the boundary values of the first-order interior vertical velocity ${u}^{1,\mathrm{I}}_3$ at the bottom boundary ($z=B$) and top boundary ($z=B+2$):
\begin{align}
    &\bm{u}^{1,\mathrm{I}}
    \cdot\bm{n}|_{z=B} =
    \cos\gamma	
	\int_{0}^{+\infty}
	\tilde{\nabla}_{\delta^0}
	\cdot\boldsymbol{u}^{0,\mathrm{BL,b}}
	d\tilde{z}'
    \label{D-3.20}\\
    &\bm{u}^{1,\mathrm{I}}
    \cdot\bm{n}|_{z=B+2} =-\cos\gamma\int_{0}^{+\infty}
	\bar{\nabla}_{\delta^0}
	\cdot\boldsymbol{u}^{0,\mathrm{BL,t}}
	d\bar{z}'. \label{D-3.21}
\end{align}

Next, applying the horizontal curl operator $\nabla_h^\perp \cdot$ to equation \eqref{D-3.25} and using the $\mathcal{O}(\delta)$ order interior incompressibility condition
$$\nabla_h \cdot (\varepsilon^{-1}\delta \boldsymbol{u}_h^{1,\mathrm{I}}) = \nabla_h \cdot {\bar{\boldsymbol{u}}}_h^1 = -\varepsilon^{-1}\delta \partial_z u_3^{1,\mathrm{I}},$$
we obtain
\begin{equation}\label{D-3.26}
    \partial_t (\nabla_h^\perp \cdot {\boldsymbol{u}}^{0,\mathrm{I}}_h) + ({\boldsymbol{u}}^{0,\mathrm{I}}_h \cdot \nabla_h)(\nabla_h^\perp \cdot {\boldsymbol{u}}^{0,\mathrm{I}}_h) = \varepsilon^{-1}\delta \partial_z u_3^{1,\mathrm{I}}.
\end{equation}
Subsequently, integrating equation \eqref{D-3.26} over $z \in [B, B+2]$, the right-hand side represents the Ekman suction effect:
\begin{equation*}
    \int_B^{B+2} \varepsilon^{-1}\delta \partial_z u_3^{1,\mathrm{I}} \,dz = \varepsilon^{-1}\delta \big( u_3^{1,\mathrm{I}}|_{z=B+2} - u_3^{1,\mathrm{I}}|_{z=B} \big).
\end{equation*}
Substituting the previously derived Ekman suction matching conditions \eqref{D-3.20}–\eqref{D-3.21} into equation \eqref{D-3.26} yields the following vorticity equation:
\begin{equation*}
    \nabla_h^\perp \cdot \bigg[  \partial_t {\boldsymbol{u}}^{0,\mathrm{I}}_h + ({\boldsymbol{u}}^{0,\mathrm{I}}_h \cdot \nabla_h) {\boldsymbol{u}}^{0,\mathrm{I}}_h + \sqrt{\tfrac{\nu}{2}} \cos^{-\frac{3}{2}}\gamma (\cos\gamma \boldsymbol{H}_0  - \boldsymbol{E}_1 ){\boldsymbol{u}}^{0,\mathrm{I}}_h  + \nabla_h B \,D_t {u}^{0,\mathrm{I}}_3\bigg] = 0,
\end{equation*}
The correspondence between this vorticity equation and the first-order velocity components in \eqref{D-3.24} is as follows: the generalized Ekman suction term $\nabla_h^\perp \cdot \big(\sqrt{\tfrac{\nu}{2}} \cos^{-\frac{3}{2}}\gamma (\cos\gamma \boldsymbol{H}_0  - \boldsymbol{E}_1 ){\boldsymbol{u}}^{0,\mathrm{I}}_h\big)$ originates from the 2D component ${\bar{\boldsymbol{u}}}_h^{1}(t,x,y)$; while the topographic feedback term $\nabla_h^\perp \cdot \big(\nabla_h B \,D_t {u}^{0,\mathrm{I}}_3\big)$ is derived from the $z$-dependent shear term $- z \nabla_h^\perp \big(D_t {u}^{0,\mathrm{I}}_3 \big)$.

Thus, we obtain the expressions for the first-order interior velocities:
\begin{align}
    \boldsymbol{u}_h^{1,\mathrm{I}} =&
\tfrac{1}{\sqrt{2}}(\cos\gamma \boldsymbol{E}_1 \boldsymbol{H}_0   + \boldsymbol{E}) {\boldsymbol{u}}^{0,\mathrm{I}}_h
-\varepsilon \delta^{-1}   
\big(
    \nabla_h^\perp B  \,D_t {u}^{0,\mathrm{I}}_3 + z\nabla_h^\perp \big( D_t {u}^{0,\mathrm{I}}_3\big)
\big)
,\label{D-3.28}\\
    {u}_3^{1,\mathrm{I}} =&
    \nabla_h B\cdot\boldsymbol{u}_h^{1,\mathrm{I}}+\tfrac{\varepsilon
    (B+1-z)}{\delta}
\nabla_h^\perp\cdot
\big(\sqrt{\tfrac{\nu}{2}} (\cos\gamma)^{-3/2} (\cos\gamma  \boldsymbol{H}_0   - \boldsymbol{E}_1){\boldsymbol{u}}^{0,\mathrm{I}}_h
    \big).
    \label{D-3.29}
\end{align}

In summary, we finally derive the 2D limiting system:
\begin{equation}\label{D-3.27}
    \begin{cases}
        \partial_t {\boldsymbol{u}}^{0,\mathrm{I}}_h + ({\boldsymbol{u}}^{0,\mathrm{I}}_h \cdot \nabla_h) {\boldsymbol{u}}^{0,\mathrm{I}}_h + \nabla_h {\bar{p}}
        + \nabla_h B \,  D_t {u}^{0,\mathrm{I}}_3 \\\qquad\quad+ \sqrt{\frac{\nu}{2}} \cos^{-\frac{3}{2}}\gamma (\cos\gamma \boldsymbol{H}_0  - \boldsymbol{E}_1 ){\boldsymbol{u}}^{0,\mathrm{I}}_h   = \boldsymbol{0}, \\
        \nabla_h \cdot {\boldsymbol{u}}^{0,\mathrm{I}}_h = 0.
    \end{cases}
\end{equation}

\begin{remark}\label{D-rem2}
    Unlike the zero-order vertical velocity ($u_3^{0,\mathrm{I}} = \nabla_h B \cdot \boldsymbol{u}_h^{0,\mathrm{I}}$), which is solely governed by topographic effects, the intrinsic mechanism of the first-order interior vertical component $u_3^{1,\mathrm{I}}$ is more complex: it includes not only the topographic effect induced by the first-order horizontal flow field ($\nabla_h B \cdot \boldsymbol{u}_h^{1,\mathrm{I}}$), but also the Ekman suction component induced by viscous friction within the boundary layer $\nabla_h^\perp\cdot\big(\sqrt{\frac{\nu}{2}} (\cos\gamma)^{-3/2} (\cos\gamma  \boldsymbol{H}_0   - \boldsymbol{E}_1){\boldsymbol{u}}^{0,\mathrm{I}}_h \big)$.
\end{remark}

\subsubsection{The $\mathcal{O}(1)$ Order Boundary Terms in the Asymptotic Expansion}\label{D-subsubsec3.2.4}

To solve the $\mathcal{O}(\delta)$ order boundary layer correction terms, we perform diagonalization on the $\mathcal{O}(1)$ order boundary part of the approximate equations. Taking the horizontal velocity $\boldsymbol{u}_h^{1,\mathrm{BL,b}}$ in the bottom boundary layer as an example, it satisfies the following inhomogeneous equations:
\begin{align}\label{D-3.30}
    \partial_{\tilde z}^2 \boldsymbol u_h^{1,\mathrm{BL,b}}
    &+
    \boldsymbol{Q}\begin{pmatrix}
        i&0\\
        0&-i
    \end{pmatrix}\boldsymbol{Q}^{-1}\boldsymbol u_h^{1,\mathrm{BL,b}}\\
    &=\boldsymbol{H_0}^{-1}
    \Bigg(\begin{pmatrix}
        1&0&B_x\\
        0&1&B_y 
    \end{pmatrix}\boldsymbol{D}^{1,\mathrm{BL,b}}+\partial_{\tilde{z}} D^{1,\mathrm{BL,b}}_0\nabla_h B\Bigg)=:
    \begin{pmatrix}
        G^{1,\mathrm{BL,b}}_1\\ 
        G^{1,\mathrm{BL,b}}_2
    \end{pmatrix},\notag
\end{align}
where the diagonalization similarity transformation matrix $\boldsymbol{Q}$ and its inverse $\boldsymbol{Q}^{-1}$ are given by:
\begin{equation*}
    \boldsymbol{Q}=\begin{pmatrix}
        \cos\gamma(1+B_y^2)
        &\cos\gamma(1+B_y^2)\\
        i-\cos\gamma B_xB_y
        &-i-\cos\gamma B_xB_y
    \end{pmatrix},
\end{equation*}
\begin{equation*}
    \boldsymbol{Q}^{-1}
    =-i
    \big(2\cos\gamma(1+B_y^2)\big)^{-1}
    \begin{pmatrix}
        i+\cos\gamma B_xB_y
        &\cos\gamma(1+B_y^2)\\
        i-\cos\gamma B_xB_y
        &-\cos\gamma(1+B_y^2)
    \end{pmatrix}.
\end{equation*}
The inhomogeneous terms on the right-hand side of system \eqref{D-3.30} are coupled from the zero-order flow field and nonlinear advection terms, where $\boldsymbol{D}^{1,\mathrm{BL,b}}$ and $D^{1,\mathrm{BL,b}}_0$ are defined respectively as:
\begin{align*}
    \boldsymbol{D}^{1,\mathrm{BL,b}}
    =&(\cos\gamma\delta)^{-1}\tilde{\nabla}_{\delta^0}(\delta p^{1,\mathrm{BL,b}})
    -\cos^2\gamma\tilde{\Delta}_{\delta^{-1}}\boldsymbol u^{0,\mathrm{BL,b}}
    \\
    &
    +\sqrt{\tfrac{\cos\gamma}{\nu}}
    \Big(
    \partial_t\boldsymbol u^{0,\mathrm{BL,b}}+
    (\boldsymbol{u}^{0}\cdot\tilde{\nabla}_{\delta^0})\boldsymbol{u}^{0,\mathrm{BL,b}}
    \Big),\\
    &
    +\sqrt{\tfrac{\cos\gamma}{\nu}}
    \Big(
    (\boldsymbol{u}_h^{0,\mathrm{BL,b}}\cdot\nabla_h)\boldsymbol{u}^{0,\mathrm{I}}+
    (\boldsymbol{u}^{1}\cdot\tilde{\nabla}_{\delta^{-1}})
    \boldsymbol{u}^{0,\mathrm{BL,b}}
    \Big),\\
    D^{1,\mathrm{BL,b}}_0=
    &\tilde{\nabla}_{\delta^0}\cdot\boldsymbol u^{0,\mathrm{BL,b}}.
\end{align*}
It should be noted that the nonlinear advection component involving $\boldsymbol{u}^{1}$ in $\boldsymbol{D}^{1,\mathrm{BL,b}}$ can be explicitly expressed using equations \eqref{D-3.18}–\eqref{D-3.19} and \eqref{D-3.28}–\eqref{D-3.29}.

Combined with the no-slip boundary conditions and the far-field decay conditions \eqref{D-3.3}–\eqref{D-3.4}, solving the inhomogeneous diagonalized system \eqref{D-3.30} yields the expression for the first-order horizontal velocity field in the bottom boundary layer:
\begin{equation}\label{D-3.31}
    \boldsymbol{u}_h^{1,\mathrm{BL,b}} = \boldsymbol{M}(\tfrac{\tilde{z}}{\sqrt{2}})
    \boldsymbol{u}_h^{1,\mathrm{I}}\big|_{z=B} 
    + \operatorname{Re}\, \boldsymbol{Q} \int_0^{\tilde{z}}\boldsymbol{Q}^{-1}
    \begin{pmatrix} 
        \mathrm{e}^{-\sqrt{i}(\tilde{z}-\tau)} {G}^{1,\mathrm{BL,b}}_1(\tau)  \\  
        \mathrm{e}^{-\sqrt{-i}(\tilde{z}-\tau)} {G}^{1,\mathrm{BL,b}}_2(\tau)  
    \end{pmatrix}\,d\tau.
\end{equation}
Similarly, the first-order horizontal correction velocity in the top boundary layer is given by:
\begin{equation}\label{D-3.32}
    \boldsymbol u_h^{1,\mathrm{BL,t}}=
    \boldsymbol{M}(\tfrac{\bar z}{\sqrt 2})
    \boldsymbol u_h^{1,\mathrm{I}}|_{z=B+2} 
    +
    \operatorname{Re}\boldsymbol{Q}
    \int_0^{\bar{z}}
    \boldsymbol{Q}^{-1}\begin{pmatrix}
        {\mathrm{e}}^{-\sqrt{-i}(\bar{z}-\tau)} G^{1,\mathrm{BL,t}}_1(\tau)
        \\[1ex]
        {\mathrm{e}}^{-\sqrt{i}(\bar{z}-\tau)} G^{1,\mathrm{BL,t}}_2(\tau) 
    \end{pmatrix}{d}\tau,
\end{equation}
where the inhomogeneous term $\boldsymbol{G}_h^{i,\mathrm{BL,t}}$ for the top boundary layer has a similar form to that for the bottom boundary layer, and can be obtained by replacing the differential operators and variables of the bottom boundary with those of the top boundary, which will not be repeated here.

Furthermore, integrating the normal momentum equations of the first-order boundary layers with the far-field conditions yields the second-order pressure terms for the top and bottom boundary layers:
\begin{align}
    p^{2,\mathrm{BL,b}}&=\int_{\tilde{z}}^{+\infty}\big(
    -\cos\gamma\partial_{\tilde{z}'}^2 u_3^{1,\mathrm{BL,b}}+\cos\gamma D_3^{1,\mathrm{BL,b}}
    \big)d\tilde{z}',\label{D-3.33}\\
    p^{2,\mathrm{BL,t}}&=\int_{\bar{z}}^{+\infty}\big(
    \cos\gamma\partial_{\bar{z}'}^2 u_3^{1,\mathrm{BL,t}}-\cos\gamma D_3^{1,\mathrm{BL,t}}
    \big)d\bar{z}'.\label{D-3.34}
\end{align}

\subsection{Divergence Correction of Approximate Solutions}\label{D-subsec3.3}

Recall that the multi-scale approximate solutions constructed in Section \eqref{D-subsec3.2} satisfy the no-slip boundary conditions:
\begin{equation}\label{D-3.35}
	\Big(\sum_{i=0}^{1}\delta^i(\boldsymbol{u}^{i,\mathrm{I}} + (1-\chi)\boldsymbol{u}^{i,\mathrm{BL,b}} + \chi\boldsymbol{u}^{i,\mathrm{BL,t}})\Big)|_{\partial\Omega}=\boldsymbol{0}.
\end{equation}
However, these approximate solutions do not satisfy the incompressibility condition, as the local divergences of some higher-order terms are not canceled during the boundary layer construction and matching process:
\begin{align}\label{D-3.36}
    &\nabla\cdot \Big(\sum_{i=0}^{1}\delta^i(\boldsymbol{u}^{i,\mathrm{I}} + (1-\chi)\boldsymbol{u}^{i,\mathrm{BL,b}} + \chi\boldsymbol{u}^{i,\mathrm{BL,t}})\Big)\\
    &= \chi \bar{\nabla}_{\delta^0}\cdot(\delta \boldsymbol{u}^{1,\mathrm{BL,t}})+(1-\chi)\tilde{\nabla}_{\delta^0}\cdot(\delta \boldsymbol{u}^{1,\mathrm{BL,b}})\neq 0.\notag
\end{align}
Moreover, cross terms involving the derivative of the cutoff function $\chi'$ arise when computing the divergence. Since $\chi'$ is non-zero only in the transition region away from the walls, where the boundary layer profiles decay exponentially, the coupling terms between $\chi'$ and the boundary layers vanish in the asymptotic expansion.

Since the right-hand side of equation \eqref{D-3.36} is of order $\mathcal{O}(\delta)$, we introduce a first-order velocity correction term $\boldsymbol{u}^c$ satisfying the following constraints:
\begin{equation}\label{D-3.37}
    \begin{cases}
        \nabla \cdot \boldsymbol{u}^{c} = \chi \bar{\nabla}_{\delta^0}\cdot(\delta \boldsymbol{u}^{1,\mathrm{BL,t}})+(1-\chi)\tilde{\nabla}_{\delta^0}\cdot(\delta \boldsymbol{u}^{1,\mathrm{BL,b}})=:f^{c}, \\
        \boldsymbol{u}^{c}|_{\partial\Omega} = \boldsymbol{0}.
    \end{cases}
\end{equation}
The following proposition holds for the divergence correction field $\boldsymbol{u}^c$:
\begin{proposition}\label{D-prop1}
    Assume that the boundary function $B(x,y)$ of the domain $\Omega$ satisfies $B(x,y)\in W^{4,\infty}(\mathbb{R}^2)$, and the right-hand side of equation \eqref{D-3.37} satisfies $f^{c} \in W^{1,\infty}(\mathbb{R}_+; L^2(\Omega))$ $(i=0,1)$. Then there exists a vector field $\boldsymbol{u}^{c} \in W^{1,\infty}(\mathbb{R}_+; H^1(\Omega))$ satisfying constraints \eqref{D-3.37}, and the following estimate holds:
    \begin{equation}\label{D-3.38}
        \lVert\boldsymbol{u}^{c}\rVert_{W^{1,\infty}(\mathbb{R}_+;H^1(\Omega))} \lesssim \lVert f^{c} \rVert_{W^{1,\infty}(\mathbb{R}_+;L^2(\Omega))}.
    \end{equation}
\end{proposition}
\begin{proof}
    According to the results on the divergence equation $\nabla \cdot \boldsymbol{u}^{c} = f^{c}$ in Refs. \refcite{Galdi1994} and \refcite{Ukai1986}, if $f^{c} \in L^2(\Omega)$, then there exists a vector field $\boldsymbol{u}^{c} \in H^1_0(\Omega)$ satisfying the homogeneous Dirichlet boundary condition, and $\lVert\boldsymbol{u}^{c}\rVert_{H^1} \lesssim \lVert f^{c} \rVert_{L^2}$. This regularity result can be extended to the $W^{1,\infty}$ space by differentiating with respect to time. Since establishing the $L^2$ norm of $f^{c}$ requires more refined boundary layer estimates, the detailed proof will be presented in subsequent subsections.
\end{proof}

Finally, we construct the corresponding first-order pressure correction term $p^{c}$ for the introduced velocity correction term. Specifically, there exists $\nabla p^{c} \in H_0^1(\Omega)$ satisfying:
\begin{align}
    \varepsilon^{-1}\nabla  p^{c} &= \nu\varepsilon \Delta\boldsymbol{u}^{c} -\boldsymbol{R} \boldsymbol{u}^{c}     - {\boldsymbol{u}}^{0}\cdot\nabla\boldsymbol{u}^{c}   \label{D-3.40}\\
    &\quad - \delta^{-1} \big( (1-\chi)\boldsymbol{u}^{c}\cdot\tilde{\nabla}_{\delta^{-1}}\boldsymbol{u}^{0,\mathrm{BL,b}} + \chi\boldsymbol{u}^{c}\cdot\bar{\nabla}_{\delta^{-1}}\boldsymbol{u}^{0,\mathrm{BL,t}} \big). \notag
\end{align}
This completes the correction of the incompressibility condition for the system.

\subsection{Estimates of Approximate Solutions}\label{D-subsec3.4}

In summary, from the asymptotic analysis of each order of the multi-scale expansion and the divergence correction, the form of the approximate solution $\big(\boldsymbol{u}_\mathrm{app}^\varepsilon,{p}_\mathrm{app}^\varepsilon\big)$ is finally established as:
\begin{equation}\label{D-3.41}
\begin{cases}
    \boldsymbol{u}_{\text{app}}^\varepsilon  = \sum_{i=0}^1\delta^i
    \big(
    \boldsymbol{u}^{i,\mathrm{I}} + (1-\chi)\boldsymbol{u}^{i,\text{BL,b}} +\chi \boldsymbol{u}^{i,\text{BL,t}} 
    \big)+\boldsymbol{u}^{c}, \\
    p_{\text{app}}^\varepsilon  = \sum_{i=0}^1 
    \Big(
        \delta^i p^{i,\mathrm{I}}
        +
        \delta^{i+1}\big((1-\chi) p^{i+1,\text{BL,b}} +\chi p^{i+1,\text{BL,t}} \big)
    \Big)+p^{c} ,
\end{cases}
\end{equation}
where the zero-order interior field $(\boldsymbol{u}^{0,\mathrm{I}},p^{1,\mathrm{I}})$ in the approximate solution \eqref{D-3.41} is determined by equations \eqref{D-3.12}, \eqref{D-3.22} and \eqref{D-3.27}, and the remaining interior terms are given by equations \eqref{D-3.8} and \eqref{D-3.28}–\eqref{D-3.29}. The Ekman correction profiles of each order boundary layer correspond to equations \eqref{D-3.10}–\eqref{D-3.11}, \eqref{D-3.13}–\eqref{D-3.16}, \eqref{D-3.18}–\eqref{D-3.19} and \eqref{D-3.31}–\eqref{D-3.34}. In addition, the divergence and pressure correction terms are given by equations \eqref{D-3.37} and \eqref{D-3.40}.

Substituting the approximate solution \eqref{D-3.41} into the Navier-Stokes-Coriolis system yields the following approximate system with a high-order remainder $\boldsymbol \rho^\varepsilon$:
\begin{equation}\label{D-3.42}
	\begin{cases}
		\begin{array}{lr}
			\partial_{t} {\boldsymbol u}_{app}^{\varepsilon}-\nu\varepsilon\Delta{\boldsymbol u}_{app}^{\varepsilon}+({\boldsymbol u}_{app}^{\varepsilon} \cdot \nabla) {\boldsymbol u}_{app}^{\varepsilon}+\varepsilon^{-1}\boldsymbol{R} {\boldsymbol u}_{app}^{\varepsilon}
			+\varepsilon^{-1} \nabla{p}_{app}^{\varepsilon}=\boldsymbol \rho^\varepsilon\,,\\
			\nabla\cdot {\boldsymbol u}_{app}^{\varepsilon}=0\,,\quad
			\boldsymbol{u}^\varepsilon_{app}|_{\partial\Omega}=0\,,
		\end{array}
	\end{cases}
\end{equation}
where the remainder $\boldsymbol{\rho}^\varepsilon=\boldsymbol \rho^{\varepsilon,1}+\boldsymbol \rho^{\varepsilon,2}$ contains the uncanceled high-order small quantities and nonlinear cross terms from the multi-scale expansion, and its explicit expansion is given by:
\begin{align}
    \boldsymbol{\rho}^{\varepsilon,1}=
    &
    \partial_t\big(\delta\boldsymbol{u}^{1}+\boldsymbol{u}^{c}\big)
    -\nu\varepsilon
    \big(
    (1-\chi)\tilde{\Delta}_{\delta^{-1}}\boldsymbol u^{1,\mathrm{BL,b}}
    +\chi\bar{\Delta}_{\delta^{-1}}\boldsymbol u^{1,\mathrm{BL,t}}
    \big)\label{D-3.43}
    \\
    &
    -\nu\varepsilon 
    \sum_{i=0}^{1}
    \big(
        \Delta(\delta^i \boldsymbol{u}^{i,\mathrm{I}})
        +
        (1-\chi)\tilde{\Delta}_{\delta^{0}}
        (\delta^i \boldsymbol{u}^{i,\mathrm{BL,b}})
        +
        \chi\bar{\Delta}_{\delta^{0}}
        (\delta^i \boldsymbol{u}^{i,\mathrm{BL,t}})
    \big)
    \notag\\        
    &
    +
    {\varepsilon}^{-1}\big(
	(1-\chi)\tilde{\nabla}_{\delta^0}
	(\delta^2p^{2,\mathrm{BL,b}})
	+\chi\bar{\nabla}_{\delta^0}(\delta^2p^{2,\mathrm{BL,t}})
	\big)
    \,,\notag
    \\[6pt]
    \boldsymbol{\rho}^{\varepsilon,2}=
	&(\delta\boldsymbol{u}^{1}+\boldsymbol{u}^{c})
	\cdot
	\big(
	\chi\bar{\nabla}_{\delta^{-1}}\boldsymbol{u}^{1,\mathrm{BL,t}}
	+(1-\chi)
	\tilde{\nabla}_{\delta^{-1}}\boldsymbol{u}^{1,\mathrm{BL,b}}
	+
    \nabla \boldsymbol{u}^{c}\big)
	\label{D-3.44}\\[2pt]
	&+
    \sum_{1\leqslant i+j\leqslant 2} 
	\delta^i\boldsymbol{u}^{i}\cdot\Big(
	\nabla (\delta^i\boldsymbol{u}^{j,\mathrm{I}})
	+
	\chi\bar{\nabla}_{\delta^{0}}(\delta^j\boldsymbol{u}^{j,\mathrm{BL,t}})
	+(1-\chi)
	\tilde{\nabla}_{\delta^{0}}(\delta^j\boldsymbol{u}^{j,\mathrm{BL,b}})
	\Big)
    \notag\\
    &+
    \sum_{i=1}^2 
    \boldsymbol{u}^{c}\cdot\Big(
	\nabla (\delta^i\boldsymbol{u}^{j,\mathrm{I}})
	+
	\chi\bar{\nabla}_{\delta^{0}}(\delta^j\boldsymbol{u}^{j,\mathrm{BL,t}})
	+(1-\chi)
	\tilde{\nabla}_{\delta^{0}}(\delta^j\boldsymbol{u}^{j,\mathrm{BL,b}})
	\Big).
	\notag
\end{align}

From the construction of the approximate solution \eqref{D-3.41}, all terms are composed of the interior flow field $\boldsymbol{u}^{0,\mathrm{I}}$. For notational simplicity, we denote $\boldsymbol{u}^{0,\mathrm{I}} = \bar{\boldsymbol{u}}$. This macroscopic limit flow field satisfies the following 2D limiting dynamical system with anisotropic geometric damping:
\begin{equation}\label{D-3.45}
\begin{cases}
    \partial_{t} \bar {\boldsymbol u}_h+(\bar {\boldsymbol u}_h \cdot \nabla_h) \bar {\boldsymbol u}_h + \nabla_h B \, D_t \bar{u}_3  +\sqrt{\frac{\nu}{2\cos\gamma}} (\boldsymbol{H}_0 -\cos^{-1}\gamma\boldsymbol{E}_1) \bar{\boldsymbol u}_h+\nabla_h\bar p = \boldsymbol{0},\\
    \bar u_3=\nabla_h{B}\cdot\bar{\boldsymbol{u}}_h,\\
    \nabla_h \cdot \bar{\boldsymbol u}_h=0.
\end{cases}
\end{equation}

This subsection aims to establish the a priori norm estimates for each term of the approximate solution $\boldsymbol{u}_\mathrm{app}^\varepsilon$ in the $L^2$ framework, and thus complete the proof of Theorem \ref{D-th1}. The specific steps are as follows: first, we present the a priori estimates for the macroscopic limit equation \eqref{D-3.45}; second, based on the regularity of the limit flow field $\bar{\boldsymbol{u}}$, we derive the norm estimates for each term of the approximate solution in the corresponding Sobolev spaces sequentially; finally, we give the $L^2$ upper bound for the remainder $\boldsymbol{\rho}^\varepsilon$ of the approximate system.

\subsubsection{A Priori Estimates for the Limiting System \eqref{D-3.45}}\label{D-subsubsec3.4.1}

To carry out the subsequent $L^2$ error analysis, we state here the a priori estimates for the limiting system \eqref{D-3.45} (detailed proofs are given in the Appendix). We first establish the decay estimates for the velocity field $\bar{\boldsymbol{u}}_{h}$, the macroscopic vorticity $\bar{\omega}:= \nabla_h^\perp\cdot\bar{\boldsymbol{u}}_h$, and the generalized vorticity $\tilde{\omega} := \nabla_h^\perp \cdot (\boldsymbol{H}_0 \bar{\boldsymbol{u}}_h)$.

\begin{proposition}\label{D-prop2}
Let the initial value $\bar{\boldsymbol{u}}_{0,h}\in H^{1}(\mathbb{R}^2)$ be a divergence-free field, with initial vorticity $\bar\omega_0=\nabla_h^\perp\cdot\bar{\boldsymbol{u}}_{0,h}\in L^{2}(\mathbb{R}^2)$. Assume that $(\bar{\boldsymbol{u}},\bar{p})$ is a solution to the limiting system \eqref{D-3.45} with initial value $\bar{\boldsymbol{u}}_0=(\bar{\boldsymbol{u}}_{0,h},\nabla_h B\cdot\bar{\boldsymbol{u}}_{0,h})$. Assume further that the smooth and globally bounded boundary surface satisfies the following geometric constraints:
\begin{align}
    &0< C_1\leqslant \cos\gamma \leqslant C_2<1,\label{D-3.46}\\
    &\sup_{(x,y)\in\mathbb{R}^2} \max\big(|k_1(x,y)|, |k_2(x,y)|\big) \leqslant C_0,\label{D-3.47}
\end{align}
where $C_0,C_1,C_2< \infty$ are positive constants, $\gamma$ is the angle between the boundary normal vector and the vertical direction, and $k_1, k_2$ are the principal curvatures of the boundary surface. Then the following decay estimates hold for the limiting system:
\begin{align*}
\|\bar{\boldsymbol u}(t)\|_{L^2(\mathbb{R}^2)}^2
&\lesssim
\|\bar{\boldsymbol{u}}_0\|_{L^2(\mathbb{R}^2)}^2  \mathrm{e}^{-{\sqrt{2\nu}}t},\\
\|\bar{\omega}(t)\|^2_{L^2(\mathbb{R}^2)} &\lesssim \Big(\|\bar{\omega}_0\|^2_{L^2(\mathbb{R}^2)}+\|\bar{\boldsymbol{u}}_{0,h}\|^2_{L^2(\mathbb{R}^2)}   \Big) \mathrm{e}^{-\sqrt{\frac{\nu}{2}}t}.
\end{align*}
\end{proposition}

\begin{proposition}\label{D-prop3}
Assume that the boundary surface $B(x,y)\in W^{2,\infty}(\mathbb{R}^2)$ and the initial flow field satisfies $\bar{\boldsymbol{u}}_{0,h}\in W^{1,\infty}(\mathbb{R}^2)$. Define
\begin{equation*}
    \mathcal{I}_0 := \max \Big\{ \|\bar{\boldsymbol{u}}_{0,h}\|_{L^2(\mathbb{R}^2)}, \, \|\bar{\boldsymbol{u}}_{0,h}\|_{L^\infty(\mathbb{R}^2)}, \, \|\bar{\omega}_0\|_{L^2(\mathbb{R}^2)}, \, \|\bar{\omega}_0\|_{L^\infty(\mathbb{R}^2)} \Big\}.
\end{equation*}
Then the following $L^\infty$ decay estimates hold for the limiting system \eqref{D-3.45}:
\begin{equation*}
    \|\bar{\boldsymbol{u}}_h(t)\|_{L^\infty(\mathbb{R}^2)} , \|\tilde{\omega}(t)\|_{L^\infty(\mathbb{R}^2)} \lesssim \mathcal{I}_0 \, \mathrm{e}^{-\frac12\sqrt{\frac{\nu}{2}}t}.
\end{equation*}
\end{proposition}

For the purpose of error analysis in the asymptotic expansion, the limit velocity field $\bar{\boldsymbol{u}}_h$ must have at least the high-order regularity of $H^3(\mathbb{R}^2)$, which leads to the following proposition.

\begin{proposition}\label{D-prop4}
If the initial value satisfies $\bar{\boldsymbol{u}}_{0,h} \in H^3(\mathbb{R}^2)$ and the boundary function $B(x,y)\in W^{4,\infty}(\mathbb{R}^2)$, then the high-order Sobolev energy of the limit velocity field and the $L^\infty$ norm of its spatial gradient satisfy the following decay estimates:
\begin{equation*}
    \|\bar{\boldsymbol{u}}_h(t)\|_{H^3(\mathbb{R}^2)} , \|\nabla_h \bar{\boldsymbol{u}}_h(t)\|_{L^\infty(\mathbb{R}^2)} \lesssim \mathrm{e}^{-\frac{1}{2}\sqrt{\frac{\nu}{2}} t}.
\end{equation*}
\end{proposition}

\subsubsection{Estimates for the Approximate Solution and the Remainder of the Approximate System}\label{D-subsubsec3.4.2}
Based on the analysis of the properties of the limit solution in the previous subsection, this subsection completes the norm estimates for the remaining expansion terms of the approximate solution and the control of the remainder of the approximate system.

\begin{proposition}[Estimates for the Zero-Order Boundary Layer Terms and Correction Terms]\label{D-prop5}
Assume that the initial flow field satisfies $\bar{\boldsymbol{u}}_{0,h} \in W^{1,\infty}(\mathbb{R}^2)$ and the boundary $B(x,y)\in W^{4,\infty}(\mathbb{R}^2)$. Then the zero-order boundary layer terms and the divergence correction term satisfy the following a priori estimates:
\begin{align*}
    \Big\lVert d(z)^{\frac{1}{2}} \big|\chi\boldsymbol u^{0,\mathrm{BL,t}}+(1-\chi)\boldsymbol u^{0,\mathrm{BL,b}}\big| \Big\rVert_{{L}^2([B,B+2];{L}^\infty(\mathbb{R}^2))}
    &\leqslant 2\sqrt{\nu}\varepsilon
    \cos^{-\frac52}\gamma
    \, \mathcal{I}_0,
\end{align*}
where $d(z)$ denotes the distance to the boundary $B(x,y)$.
\end{proposition}
\begin{proof}
Based on the expressions for the zero-order boundary layer terms of the top and bottom layers in equations \eqref{D-3.10}–\eqref{D-3.11}, \eqref{D-3.13} and \eqref{D-3.15}, we have the following estimates:
\begin{align*}
    |\boldsymbol{u}^{0,\mathrm{BL,b}}|^2 &= |\boldsymbol{u}_h^{0,\mathrm{BL,b}}|^2 + |u_3^{0,\mathrm{BL,b}}|^2 \\
    &\leqslant |\boldsymbol{u}_h^{0,\mathrm{BL,b}}|^2 \big( 1 + |\nabla_h B|^2 \big) \\
    &\leqslant \mathrm{e}^{-{\sqrt{2}\tilde{z}}} ( \cos^{-2}\gamma + \|\boldsymbol{H}_0\|^2 ) |\bar{\boldsymbol{u}}_h|^2\\
    &\leqslant 2\mathrm{e}^{-{\sqrt{2}\tilde{z}}}   \cos^{-2}\gamma   |\bar{\boldsymbol{u}}_h|^2
\end{align*}
and
\begin{equation*}
    |\boldsymbol{u}^{0,\mathrm{BL,t}}|^2 \leqslant 2\mathrm{e}^{-{\sqrt{2}\bar{z}}}   \cos^{-2}\gamma   |\bar{\boldsymbol{u}}_h|^2.
\end{equation*}
Furthermore, using the fast-variable decay structure of the boundary layers, which is of the form $\mathrm{e}^{-\frac{\sqrt{2}(2+B-z)}{\delta}}$ and $\mathrm{e}^{-\frac{\sqrt{2}(z-B)}{\delta}}$, we can deduce that:
\begin{align*}
    &\Big\lVert d(z)^{\frac{1}{2}} |\chi\boldsymbol u^{0,\mathrm{BL,t}}+(1-\chi)\boldsymbol u^{0,\mathrm{BL,b}}| \Big\rVert_{{L}^2([B,B+2];{L}^\infty(\mathbb{R}^2))}\\\notag
    &\leqslant
    \tfrac{\delta}{\cos\gamma}\lVert\bar{\boldsymbol{u}}_h\rVert_{{L}^\infty(\mathbb{R}^2)}
    \Big(
    \int_{[0,\frac{2\sqrt{2}}{\delta}]}
    \tfrac{\sqrt{2}(z-B)}{\delta}
    \mathrm{e}^{-\frac{\sqrt{2}(z-B)}{\delta}}
    \,{d}\tfrac{\sqrt{2}(z-B)}{\delta}    
    \Big)^{\frac12}\\
    &\quad+\tfrac{\delta}{\cos\gamma}\lVert\bar{\boldsymbol{u}}_h\rVert_{{L}^\infty(\mathbb{R}^2)}
    \Big(\int_{[0,\frac{2\sqrt{2}}{\delta}]}
    \tfrac{\sqrt{2}(2+B-z)}{\delta}
    \mathrm{e}^{-\frac{\sqrt{2}(2+B-z)}{\delta}}\,{d}\tfrac{\sqrt{2}(2+B-z)}{\delta}\Big)^{\frac12}\\
    &\leqslant 2\sqrt{\nu}\varepsilon
    \cos^{-\frac52}\gamma
    \, \mathcal{I}_0,
\end{align*}
where the last inequality follows from Proposition \ref{D-prop3}.
\end{proof}

\begin{proposition}[$L^4$ Estimates for the First-Order Approximate Solution]\label{D-prop6}
Assume that the initial flow field satisfies $\bar{\boldsymbol{u}}_{0,h}\in W^{1,\infty}(\mathbb{R}^2)\cap H^{3}(\mathbb{R}^2)$. Then the following estimates hold for the first-order approximate solution and its time derivative:
\begin{align*}
    \lVert 
    \delta
    (
    \boldsymbol{u}^{1} 
    +\boldsymbol{u}^{c}   
    )   \rVert_{{L}^4(\Omega)}
    &\lesssim
    \varepsilon
    \big( 1+\lVert\bar{\boldsymbol{u}}_{0,h}\rVert_{W^{1,\infty}(\mathbb{R}^2)}  \big) 
    \lVert\bar{\boldsymbol{u}}_{0,h}\rVert_{H^3(\mathbb{R}^2)},\\
    \lVert 
    \delta
    \partial_t(
    \boldsymbol{u}^{1} 
    +\boldsymbol{u}^{c}   
    )   \rVert_{{L}^2(\Omega)}
    &\lesssim
    \varepsilon
    \big( 1+\lVert\bar{\boldsymbol{u}}_{0,h}\rVert_{W^{1,\infty}(\mathbb{R}^2)} \big)^2
    \lVert\bar{\boldsymbol{u}}_{0,h}\rVert_{H^3(\mathbb{R}^2)}.
\end{align*}
\end{proposition}
\begin{proof}
From Propositions \ref{D-prop2}–\ref{D-prop4}, we have $\bar{\boldsymbol{u}}_{0,h}\in W^{1,\infty}(\mathbb{R}^2)\cap H^{3}(\mathbb{R}^2)$. The result follows from standard inequality scaling based on the explicit expressions for the first-order interior terms, boundary layer terms and correction terms, combined with the high-order Sobolev embedding of the limit flow field and the limit equation \eqref{D-3.45}. The details are omitted here.
\end{proof}

Based on the results of Propositions \ref{D-prop5} and \ref{D-prop6}, when the initial velocity field satisfies $\bar{\boldsymbol{u}}_{0,h}\in W^{1,\infty}(\mathbb{R}^2)\cap H^{3}(\mathbb{R}^2)$, the remainder of the approximate system satisfies:
\begin{equation}\label{D-3.48}
    \lVert \boldsymbol \rho^\varepsilon \rVert_{L^{2}(\Omega)} 
    \lesssim\varepsilon
    \big( 1+\lVert\bar{\boldsymbol{u}}_{h}\rVert_{W^{1,\infty}(\mathbb{R}^2)} \big)^2
    \lVert\bar{\boldsymbol{u}}_{h}\rVert_{H^3(\mathbb{R}^2)}.
\end{equation}

\section{$L^2$ Error Estimates}\label{D-sec4}

In this section, we prove the following result: as the parameter $\varepsilon\rightarrow 0$, the weak solution $\boldsymbol{u}^\varepsilon$ of system \eqref{D-1.5} under the initial condition \eqref{D-1.10} converges to $\bar{\boldsymbol{u}}$ in the ${L}^2(\Omega)$ norm, i.e., $\lVert\boldsymbol{u}^\varepsilon-\bar{\boldsymbol{u}}\rVert_{{L}^2(\Omega)}$ tends to zero. 
Given the structure of the approximate solution $\boldsymbol u_\mathrm{app}^{\varepsilon}$ in equation \eqref{D-3.41}, it suffices to prove that $\lVert\boldsymbol u^{\varepsilon}-\boldsymbol{u}_\mathrm{app}^{\varepsilon}\rVert_{{L}^2(\Omega)}\rightarrow 0$.

Recall that $\boldsymbol{u}^{\varepsilon}$ and $\boldsymbol u_{\mathrm{app}}^{\varepsilon}$ satisfy the following systems respectively:
\begin{displaymath}
    \begin{cases}
        \partial_t \boldsymbol{u}^\varepsilon-\nu\varepsilon \Delta \boldsymbol{u}^\varepsilon+(\boldsymbol{u}^\varepsilon \cdot \nabla)\boldsymbol{u}^\varepsilon+\varepsilon^{-1}\boldsymbol{R}\boldsymbol{u}^\varepsilon+\varepsilon^{-1}\nabla p^\varepsilon=0,\\
        \nabla\cdot \boldsymbol{u}^\varepsilon=0,\quad
        \boldsymbol{u}^\varepsilon|_{\partial\Omega}=0,\\
        \boldsymbol{u}^\varepsilon|_{t=0}=\boldsymbol{u}^\varepsilon_0, 
    \end{cases}
\end{displaymath}
and
\begin{displaymath}
    \begin{cases}
        \partial_{t} {\boldsymbol u}_\mathrm{app}^{\varepsilon}-\nu\varepsilon\Delta{\boldsymbol u}_\mathrm{app}^{\varepsilon}+({\boldsymbol u}_\mathrm{app}^{\varepsilon} \cdot \nabla) {\boldsymbol u}_\mathrm{app}^{\varepsilon}+\varepsilon^{-1}\boldsymbol{R} \boldsymbol{u}_\mathrm{app}^{\varepsilon}
        +\varepsilon^{-1}\nabla{P}_\mathrm{app}^{\varepsilon}=\boldsymbol \rho^\varepsilon,\\
        \nabla\cdot \boldsymbol{u}_\mathrm{app}^\varepsilon=0,\quad
        {\boldsymbol u}_\mathrm{app}^{\varepsilon}|_{\partial\Omega}=0,\\
        \boldsymbol{u}_\mathrm{app}^\varepsilon|_{t=0}=
        \boldsymbol{u}_{\mathrm{app},0}^\varepsilon=
        \bar{\boldsymbol u}_0
        + \big(
            \delta \boldsymbol{u}^{1,\mathrm{I}}
        +\sum_{i=0}^1\delta^i
    ((1-\chi)\boldsymbol{u}^{i,\text{BL,b}} +\chi \boldsymbol{u}^{i,\text{BL,t}} 
    +\boldsymbol{v}^i   
    )\big)|_{t=0}.
    \end{cases}
\end{displaymath}
Let $\boldsymbol{u}^\varepsilon-\boldsymbol{u}^\varepsilon_\mathrm{app}=\boldsymbol{v}^\varepsilon$ and $p^\varepsilon-p_\mathrm{app}^\varepsilon=p^\varepsilon_v$. Then $(\boldsymbol{v}^\varepsilon,p^\varepsilon_v)$ satisfies:
\begin{equation}\label{D-4.1}
    \left\{\begin{aligned}
        &\partial_{t} {\boldsymbol v}^{\varepsilon}-\nu\varepsilon\Delta{\boldsymbol v}^{\varepsilon}
        +{\boldsymbol u}^{\varepsilon} \cdot \nabla {\boldsymbol v}^{\varepsilon}
        +{\varepsilon}^{-1}\boldsymbol{R}\boldsymbol v^{\varepsilon}
        +{\varepsilon}^{-1}\nabla p^\varepsilon_v
        =-{\boldsymbol v}^{\varepsilon} \cdot \nabla\boldsymbol{u}^\varepsilon_\mathrm{app}
        -\boldsymbol \rho^\varepsilon,\\
        &\nabla\cdot {\boldsymbol v}^{\varepsilon}=0,\quad{\boldsymbol v}^{\varepsilon}|_{\partial\Omega}=0,\\
        &{\boldsymbol v}^{\varepsilon}|_{t=0}=\boldsymbol{u}_0^\varepsilon-{\boldsymbol u}_{\mathrm{app},0}^{\varepsilon}.
    \end{aligned}\right.
\end{equation}

Taking the $L^2$ inner product of system \eqref{D-4.1} with $\boldsymbol v^{\varepsilon}$ yields:
\begin{align}\label{D-4.2}
    &\frac{1}{2}\frac{d}{dt}\lVert\boldsymbol v^{\varepsilon}\rVert^2_{{L}^2(\Omega)}
    +\nu\varepsilon\lVert\nabla\boldsymbol v^{\varepsilon}\rVert^2_{{L}^2(\Omega)}\\\notag
    &\leqslant|\langle{\boldsymbol u}^{\varepsilon} \cdot \nabla {\boldsymbol v}^{\varepsilon},{\boldsymbol v}^{\varepsilon}\rangle|
    +|\langle
    {\varepsilon}^{-1}\boldsymbol{R}\boldsymbol v^{\varepsilon},{\boldsymbol v}^{\varepsilon}\rangle|
    +|\langle
    {\varepsilon}^{-1}\nabla p^\varepsilon_v,{\boldsymbol v}^{\varepsilon}\rangle|
    \\\notag
    &\quad+|\langle
    {\boldsymbol v}^{\varepsilon} \cdot \nabla\boldsymbol{u}^\varepsilon_\mathrm{app},{\boldsymbol v}^{\varepsilon}\rangle|
    +|\langle
    \boldsymbol \rho^\varepsilon   ,{\boldsymbol v}^{\varepsilon}\rangle|, 
\end{align}
where $\langle \cdot, \cdot \rangle$ denotes the $L^2(\Omega)$ inner product.

First, using the incompressibility condition of $\boldsymbol{v}^\varepsilon$ and the skew-symmetry of the operator $\boldsymbol{R}$, we have:
\begin{equation}\label{D-4.3}
    |\langle{\boldsymbol u}^{\varepsilon} \cdot \nabla {\boldsymbol v}^{\varepsilon},{\boldsymbol v}^{\varepsilon}\rangle|
    =|\langle
    {\varepsilon}^{-1}\boldsymbol{R}\boldsymbol v^{\varepsilon},{\boldsymbol v}^{\varepsilon}\rangle|
    =|\langle
    {\varepsilon}^{-1}\nabla p^\varepsilon_v,{\boldsymbol v}^{\varepsilon}\rangle|=0.   
\end{equation}

Next, for the convective coupling term $|\langle {\boldsymbol v}^{\varepsilon} \cdot \nabla\boldsymbol{u}^\varepsilon_\mathrm{app},{\boldsymbol v}^{\varepsilon}\rangle|$, we first expand it as:
\begin{align}\label{D-4.4}
    &|\langle
    {\boldsymbol v}^{\varepsilon} \cdot \nabla\boldsymbol{u}^\varepsilon_\mathrm{app},{\boldsymbol v}^{\varepsilon}\rangle|\\
    &\leqslant|\langle{\boldsymbol v}_h^{\varepsilon} \cdot \nabla_h \bar{\boldsymbol u},{\boldsymbol v}^{\varepsilon}\rangle|
    +\big|\big\langle{\boldsymbol v}^{\varepsilon} \cdot \nabla
    ((1-\chi)\boldsymbol{u}^{0,\text{BL,b}} +\chi \boldsymbol{u}^{0,\text{BL,t}} ),{\boldsymbol v}^{\varepsilon}\big\rangle\big|\notag
    \\
    &\quad+\big|\big\langle{\boldsymbol v}^{\varepsilon} \cdot \nabla
    \big(
    \delta
    (
    \boldsymbol{u}^{1} 
    +\boldsymbol{u}^{c}   
    )    
    \big),{\boldsymbol v}^{\varepsilon}\big\rangle\big|.\notag
\end{align}
Using Hölder's inequality, the first two terms on the right-hand side of equation \eqref{D-4.4} can be controlled by the following estimate:
\begin{equation}\label{D-4.5}
    |\langle{\boldsymbol v}_h^{\varepsilon} \cdot \nabla_h \bar{\boldsymbol u},{\boldsymbol v}^{\varepsilon}\rangle|
    \leqslant
        \lVert \nabla\bar{\boldsymbol u}_h \rVert_{{L}^{\infty}(\mathbb{R}^2)}
        \lVert{\boldsymbol v}^{\varepsilon}\rVert^2_{{L}^2(\Omega)}.
\end{equation}
For the convective term induced by the zero-order boundary layer, direct estimation would lead to singularities due to its fast-variable structure. Therefore, we perform integration by parts on this term using the incompressibility condition of the flow field:
\begin{align*}
    &\big|\big\langle{\boldsymbol v}^{\varepsilon} \cdot \nabla
    ((1-\chi)\boldsymbol{u}^{0,\text{BL,b}} +\chi \boldsymbol{u}^{0,\text{BL,t}} ),{\boldsymbol v}^{\varepsilon}\big\rangle\big|\\
    &=\big|\big\langle{\boldsymbol v}^{\varepsilon} \cdot \nabla{\boldsymbol v}^{\varepsilon},(1-\chi)\boldsymbol{u}^{0,\text{BL,b}} +\chi \boldsymbol{u}^{0,\text{BL,t}} \big\rangle\big|\\\notag
    &=\Big|\int_{\mathbb{R}^2\times[B,B+1]} {\boldsymbol v}^{\varepsilon} \cdot \nabla {\boldsymbol v}^{\varepsilon}\cdot
    (1-\chi)\boldsymbol u^{0,\mathrm{BL,b}} 
    \,{d}x{d}y{d}z
    \Big|\\\notag
    &\quad+\Big|\int_{\mathbb{R}^2\times[B+1,B+2]} {\boldsymbol v}^{\varepsilon} \cdot \nabla {\boldsymbol v}^{\varepsilon}\cdot
    \chi\boldsymbol u^{0,\mathrm{BL,t}}
    \,{d}x{d}y{d}z\Big|.
\end{align*}
Considering that the perturbation velocity ${\boldsymbol v}^{\varepsilon}$ satisfies the no-slip boundary condition, using the Cauchy-Schwarz inequality, we derive:
\begin{align*}
    |{\boldsymbol v}^{\varepsilon}|
    =\Big| \int_B^z \partial_\varsigma {\boldsymbol v}^{\varepsilon}(t,x,y,\varsigma)\,{d}\varsigma\Big|
    \leqslant d(z)^{\frac{1}{2}}\lVert\partial_z {\boldsymbol v}^{\varepsilon}\rVert_{{L}^2([B,B+2])},
\end{align*}
where $d(z)$ denotes the distance to the boundary $B(x,y)$. 
Combining with the weighted norm estimate established in Proposition \ref{D-prop5}, this term simplifies to:
\begin{align}\label{D-4.6}
    &\big|\big\langle{\boldsymbol v}^{\varepsilon} \cdot \nabla
    ((1-\chi)\boldsymbol{u}^{0,\text{BL,b}} +\chi \boldsymbol{u}^{0,\text{BL,t}} ),{\boldsymbol v}^{\varepsilon}\big\rangle\big|
    \\
    &\leqslant\int_\Omega
    \lVert\partial_z {\boldsymbol v}^{\varepsilon}\rVert_{{L}^2([B,B+2])}
    |\nabla {\boldsymbol v}^{\varepsilon}|
    d(z)^{\frac{1}{2}}
    \big|\chi\boldsymbol u^{0,\mathrm{BL,t}}+(1-\chi)\boldsymbol u^{0,\mathrm{BL,b}}\big|\,{d}x{d}y{d}z\notag\\
    &\leqslant \lVert\nabla{\boldsymbol v}^{\varepsilon}\rVert^2_{{L}^2(\Omega)}
    \lVert d(z)^{\frac{1}{2}} \big|\chi\boldsymbol u^{0,\mathrm{BL,t}}+(1-\chi)\boldsymbol u^{0,\mathrm{BL,b}}\big| \rVert_{{L}^2([B,B+2];{L}^\infty(\mathbb{R}^2))}\notag\\
    &\leqslant 2\sqrt{\nu}\varepsilon
    \cos^{-\frac52}\gamma
    \, \lVert\bar{\boldsymbol{u}}_h\rVert_{{L}^\infty(\mathbb{R}^2)}\lVert\nabla{\boldsymbol v}^{\varepsilon}\rVert^2_{{L}^2(\Omega)}.\notag
\end{align}
For the convective error corresponding to the first-order approximate solution, using the 3D Gagliardo-Nirenberg interpolation inequality:
\[
\|{\boldsymbol v}^{\varepsilon}\|_{L^4} \lesssim \|{\boldsymbol v}^{\varepsilon}\|_{L^2}^{1/4} \|\nabla{\boldsymbol v}^{\varepsilon}\|_{L^2}^{3/4}.
\]
Combining with the $L^4$ estimate from Proposition \ref{D-prop6}, we obtain:
\begin{align}\label{D-4.7}
    &\big|\big\langle{\boldsymbol v}^{\varepsilon} \cdot \nabla
    \big(
    \delta
    (
    \boldsymbol{u}^{1} 
    +\boldsymbol{u}^{c}   
    )   
    \big),{\boldsymbol v}^{\varepsilon}\big\rangle\big|=\big|\big\langle{\boldsymbol v}^{\varepsilon} \cdot \nabla
    {\boldsymbol v}^{\varepsilon},\delta
    (
    \boldsymbol{u}^{1} 
    +\boldsymbol{u}^{c}   
    )\big\rangle\big|\\
    &\leqslant \lVert {\boldsymbol v}^{\varepsilon} \rVert_{{L}^4(\Omega)}
    \lVert \nabla{\boldsymbol v}^{\varepsilon} \rVert_{{L}^2(\Omega)}
    \lVert 
    \delta
    (
    \boldsymbol{u}^{1} 
    +\boldsymbol{u}^{c}   
    )  \rVert_{{L}^4(\Omega)}\notag\\
    &\leqslant 
    \varepsilon
    \lVert {\boldsymbol v}^{\varepsilon} \rVert^{\frac14}_{{L}^2(\Omega)}
    \lVert \nabla{\boldsymbol v}^{\varepsilon} \rVert^{\frac74}_{{L}^2(\Omega)}
    \big( 1+\lVert\bar{\boldsymbol{u}}_{h}\rVert_{W^{1,\infty}(\mathbb{R}^2)}  \big) 
    \lVert\bar{\boldsymbol{u}}_{h}\rVert_{H^3(\mathbb{R}^2)}
    \notag\\
    &\lesssim
    \frac{\nu\varepsilon}{2}\lVert \nabla{\boldsymbol v}^{\varepsilon} \rVert^{2}_{{L}^2(\Omega)}
    +\varepsilon
    \big( 1+\lVert\bar{\boldsymbol{u}}_{h}\rVert_{W^{1,\infty}(\mathbb{R}^2)}  \big)^8
    \lVert\bar{\boldsymbol{u}}_{h}\rVert^8_{H^3(\mathbb{R}^2)}\lVert {\boldsymbol v}^{\varepsilon} \rVert^{2}_{{L}^2(\Omega)}
    .\notag
\end{align}

Combining the above estimates \eqref{D-4.4}–\eqref{D-4.7}, we get:
\begin{align}\label{D-4.8}
    |\langle
    {\boldsymbol v}^{\varepsilon} \cdot \nabla\boldsymbol{u}^\varepsilon_\mathrm{app},{\boldsymbol v}^{\varepsilon}\rangle|
    \lesssim&
    \frac{\nu\varepsilon}{2} \lVert \nabla{\boldsymbol v}^{\varepsilon} \rVert^{2}_{{L}^2(\Omega)}
    +2\sqrt{\nu}\varepsilon\, \lVert\bar{\boldsymbol{u}}_h\rVert_{{L}^\infty(\mathbb{R}^2)}\lVert \nabla{\boldsymbol v}^{\varepsilon} \rVert^{2}_{{L}^2(\Omega)}\\
    &+
    (
        \lVert\bar{\boldsymbol{u}}_{h}\rVert_{W^{1,\infty}(\mathbb{R}^2)}+
        \varepsilon
    ( 1+\lVert\bar{\boldsymbol{u}}_{h}\rVert_{W^{1,\infty}(\mathbb{R}^2)}  )^8\lVert\bar{\boldsymbol{u}}_{h}\rVert^8_{H^3(\mathbb{R}^2)}
    )
    \lVert{\boldsymbol v}^{\varepsilon}\rVert^2_{{L}^2(\Omega)}.\notag
\end{align}

Finally, we handle the remainder of the approximate system. Using the $L^2$ a priori estimate for the remainder \eqref{D-3.48}, we obtain:
\begin{align}\label{D-4.9}
    |\langle\boldsymbol \rho^\varepsilon,{\boldsymbol v}^{\varepsilon}\rangle|\leqslant&
    \lVert\boldsymbol \rho^\varepsilon\rVert_{{L}^2(\Omega)}\lVert{\boldsymbol v}^{\varepsilon}\rVert_{{L}^2(\Omega)}\\
    \lesssim&
    \varepsilon\big( 1+\lVert\bar{\boldsymbol{u}}_{h}\rVert_{W^{1,\infty}(\mathbb{R}^2)}  \big)^2
    \lVert\bar{\boldsymbol{u}}_{h}\rVert_{H^3(\mathbb{R}^2)}\lVert{\boldsymbol v}^{\varepsilon}\rVert_{{L}^2(\Omega)}.\notag
\end{align} 

Substituting the convective term estimate \eqref{D-4.8} and the remainder estimate \eqref{D-4.9} back into the energy inequality \eqref{D-4.2}, we further simplify it to:
\begin{align}\label{D-4.10}
    &\frac{d}{dt}\lVert\boldsymbol v^{\varepsilon}\rVert^2_{{L}^2(\Omega)}
    +\nu\varepsilon\lVert\nabla\boldsymbol v^{\varepsilon}\rVert^2_{{L}^2(\Omega)}\\\notag
    &\lesssim
    4\sqrt{\nu}\varepsilon
    \cos^{-\frac52}\gamma
    \, \lVert\bar{\boldsymbol{u}}_h\rVert_{{L}^\infty(\mathbb{R}^2)}\lVert\nabla{\boldsymbol v}^{\varepsilon}\rVert^2_{L^2(\Omega)}
    \\\notag
    &\quad+\varepsilon\big( 1+\lVert\bar{\boldsymbol{u}}_{h}\rVert_{W^{1,\infty}(\mathbb{R}^2)}  \big)^2
    \lVert\bar{\boldsymbol{u}}_{h}\rVert_{H^3(\mathbb{R}^2)}\lVert{\boldsymbol v}^{\varepsilon}\rVert_{{L}^2(\Omega)}\\\notag
    &\quad+
    \Big(
        \lVert\bar{\boldsymbol{u}}_{h}\rVert_{W^{1,\infty}(\mathbb{R}^2)}+
        \varepsilon
    \big( 1+\lVert\bar{\boldsymbol{u}}_{h}\rVert_{W^{1,\infty}(\mathbb{R}^2)}  \big)^8\lVert\bar{\boldsymbol{u}}_{h}\rVert^8_{H^3(\mathbb{R}^2)}
    \Big)
    \lVert{\boldsymbol v}^{\varepsilon}\rVert^2_{{L}^2(\Omega)}.\notag
\end{align}
According to the assumption \eqref{D-1.11} imposed on the limit flow field in Theorem \ref{D-th2}, the gradient error term on the right-hand side is guaranteed to be controlled by viscous dissipation, i.e.,
\begin{equation}\label{D-4.11}
    4\sqrt{\nu}\varepsilon
    \cos^{-\frac52}\gamma
    \, \lVert\bar{\boldsymbol{u}}_h\rVert_{{L}^\infty(\mathbb{R}^2)}\lVert\nabla{\boldsymbol v}^{\varepsilon}\rVert^2_{L^2(\Omega)}
    <
    \nu\varepsilon\lVert\nabla\boldsymbol v^{\varepsilon}\rVert^2_{{L}^2(\Omega)}.
\end{equation}
Furthermore, for sufficiently small parameter $\varepsilon$, substituting the decay estimates for the limit flow field established in Propositions \ref{D-prop2}–\ref{D-prop4} into \eqref{D-4.10} yields:
\begin{equation}\label{D-4.12}
    \frac{d}{dt}\lVert\boldsymbol v^{\varepsilon}\rVert^2_{{L}^2(\Omega)}
    \lesssim
    \lVert{\boldsymbol v}^{\varepsilon}\rVert^2_{{L}^2(\Omega)}\mathrm{e}^{-\frac12\sqrt{\frac{\nu}{2}}  t}
    +
    \varepsilon\mathrm{e}^{-\frac12\sqrt{\frac{\nu}{2}}  t}.
\end{equation}
Finally, applying Gronwall's lemma to the above inequality over the time interval $[0, \infty)$ and combining with the initial convergence condition \eqref{D-1.10}, we obtain the final convergence result \eqref{D-1.13}.

\section{Geometric Effects of the Limiting System Under Special Boundaries}\label{D-sec5}

This section aims to analyze the multiple mechanisms driving the vertical motion of rotating fluids over $\mathcal{O}(1)$ large-undulation terrain, and by comparing with the $\mathcal{O}(\varepsilon)$ small-amplitude case, understand the influence of local slope and geometric metric on the dynamical behavior of boundary layers.
Within the framework of multi-scale asymptotic analysis, the structure of the vertical velocity field can intuitively reflect the motion law of fluids passing over non-flat terrain.

Under the small-amplitude assumption (boundary $z = \varepsilon B$), the zero-order vertical velocity is zero, and the vertical motion of the fluid only appears at the $\mathcal{O}(\varepsilon)$ order. The structure of the interior vertical velocity $u_{3,\text{small}}^{\mathrm{I}}$ is given by:
\begin{equation}\label{D-5.1}
    u_{3,\text{small}}^{\mathrm{I}} = \varepsilon \nabla_h B \cdot \bar{\boldsymbol{u}}_h + \varepsilon \sqrt{\tfrac{\nu}{2}} (1+\varepsilon B-z) \nabla_h^\perp \cdot \bar{\boldsymbol{u}}_h.
\end{equation}
The mechanisms reflected by each term are as follows:
\begin{enumerate}
    \item \textbf{Topographic lifting:} The first term represents the vertical movement of the fluid along small topographic undulations.
    \item \textbf{Classical Ekman suction:} The second term is driven by the macroscopic vorticity of the background flow field, reflecting the classical isotropic damping dissipation, and the suction intensity is not affected by the topographic geometric metric.
\end{enumerate}

In contrast, under the large-undulation non-flat assumption (boundary $z=B$), the vertical velocity field $u_3^{\mathrm{I}}$ exhibits coupling of multiple mechanisms:
\begin{align}\label{D-5.2}
    {u}_3^{\mathrm{I}} =& {\nabla_h B\cdot\bar{\boldsymbol{u}}_h} + {\delta \nabla_h B\cdot\boldsymbol{u}_h^{1,\mathrm{I}}} \\
    &+ {\varepsilon (B+1-z) \nabla_h^\perp\cdot \Big(\sqrt{\tfrac{\nu}{2}} (\cos\gamma)^{-1/2} \boldsymbol{H}_0 \bar{\boldsymbol{u}}_h \Big)}  \notag\\
    &+ {\varepsilon (B+1-z) \nabla_h^\perp\cdot \Big(\sqrt{\tfrac{\nu}{2}} 
    (\cos\gamma)^{-3/2} \bar{\boldsymbol{u}}^\bot_h \Big)}.\notag 
\end{align}
The mechanisms reflected by each term are as follows:
\begin{enumerate}
    \item \textbf{Topographic lifting terms:} The first two terms represent the $\mathcal{O}(1)$ macroscopic lifting effect dominated by topography and its $\mathcal{O}(\delta)$ order secondary correction. The macroscopic fluid undergoes significant ascending or descending motion along the boundary.
    \item \textbf{Generalized Ekman suction:} The third term represents the dynamical mass exchange induced by bottom friction. The introduction of the metric tensor $\boldsymbol{H}_0$ and local slope $\cos\gamma$ characterizes the anisotropic dissipation caused by topographic curvature and leads to spatial inhomogeneity of the suction intensity.
    \item \textbf{Rotational deflection:} The last term indicates that the system still retains part of the Coriolis effect in the limit state, which is controlled by the topographic geometric features.
\end{enumerate}

Next, we set the boundary function as a one-dimensional ridge $B(x) = H_m \exp(-x^2/a^2)$, where $H_m$ is the peak height. This topographic gradient divides the region into a climbing zone ($x < 0$, $B'>0$) and a sinking zone ($x > 0$, $B'<0$). Meanwhile, we set the background zonal jet crossing the ridge as $\bar{\boldsymbol{u}}_h = (U_m \exp(-y^2/L_y^2), 0)^T$, where $U_m$ determines the maximum wind speed and $L_y$ characterizes the lateral coverage of the geostrophic wind. To intuitively present the evolution characteristics, we set the dimensionless viscosity $\nu = 1.0$, perturbation parameter $\varepsilon = 0.05$ in the numerical calculation, and take the spatial characteristic parameters of the jet and topography as $U_m = 1.0, L_y = 2.0$ and $H_m = 2.0, a = 2.0$.

\subsection{Comparison Between Small-Amplitude and Large-Undulation Models}

\begin{figure}[htbp]
    \centering
    \includegraphics[width=0.9\textwidth]{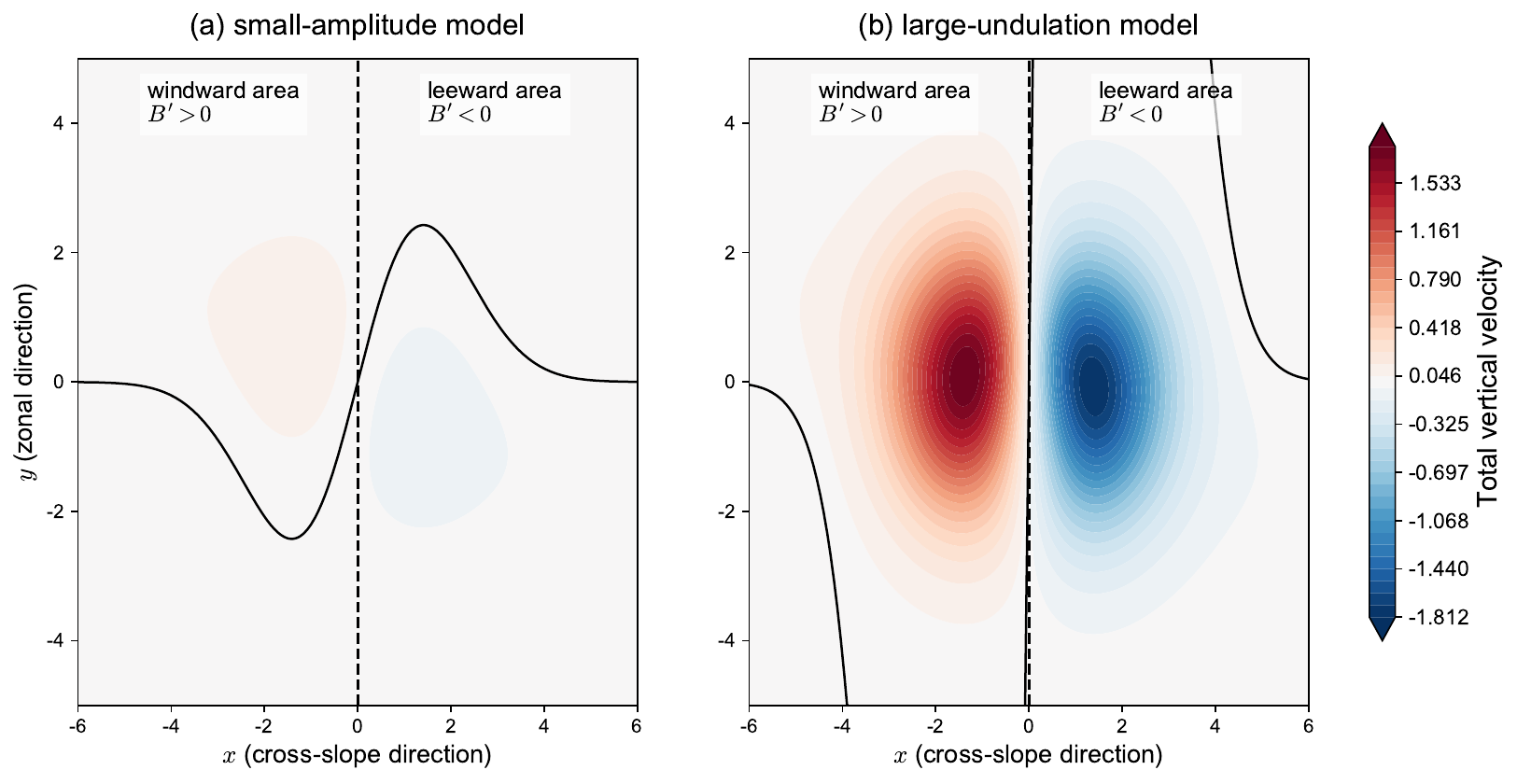}
    \caption{\textbf{Comparison of total vertical velocity fields between the small-amplitude model and the large-undulation model}}
    \label{fig:macro_compare}
\end{figure}

First, under the same zonal jet background, we compare the total vertical velocity field distributions of the $\mathcal{O}(\varepsilon)$ small-amplitude model and the $\mathcal{O}(1)$ large-undulation model in this paper (see Fig. \ref{fig:macro_compare}). The black solid lines in the figure mark the zero-contour lines of vertical velocity ($u_3=0$) to distinguish the ascending and descending regions of the flow field.
As shown in Fig. \ref{fig:macro_compare}(a), the vertical velocity of the small-amplitude model is at the $\mathcal{O}(\varepsilon)$ perturbation level. In contrast, under large-undulation terrain (Fig. \ref{fig:macro_compare}(b)), the total vertical velocity $u_3^{\mathrm{I}}$ changes significantly. Here, the $\mathcal{O}(1)$ order topographic kinematic forcing determines the ascending and descending motion of the fluid. The magnitude comparison between the two indicates that for large-scale undulating terrain, the impermeable boundary constraint plays a dominant role in vertical mass exchange.

\subsection{Superposition Process of Multiple Dynamical Mechanisms in the Large-Undulation Model}

\begin{figure}[htbp]
    \centering
    \includegraphics[width=1\textwidth]{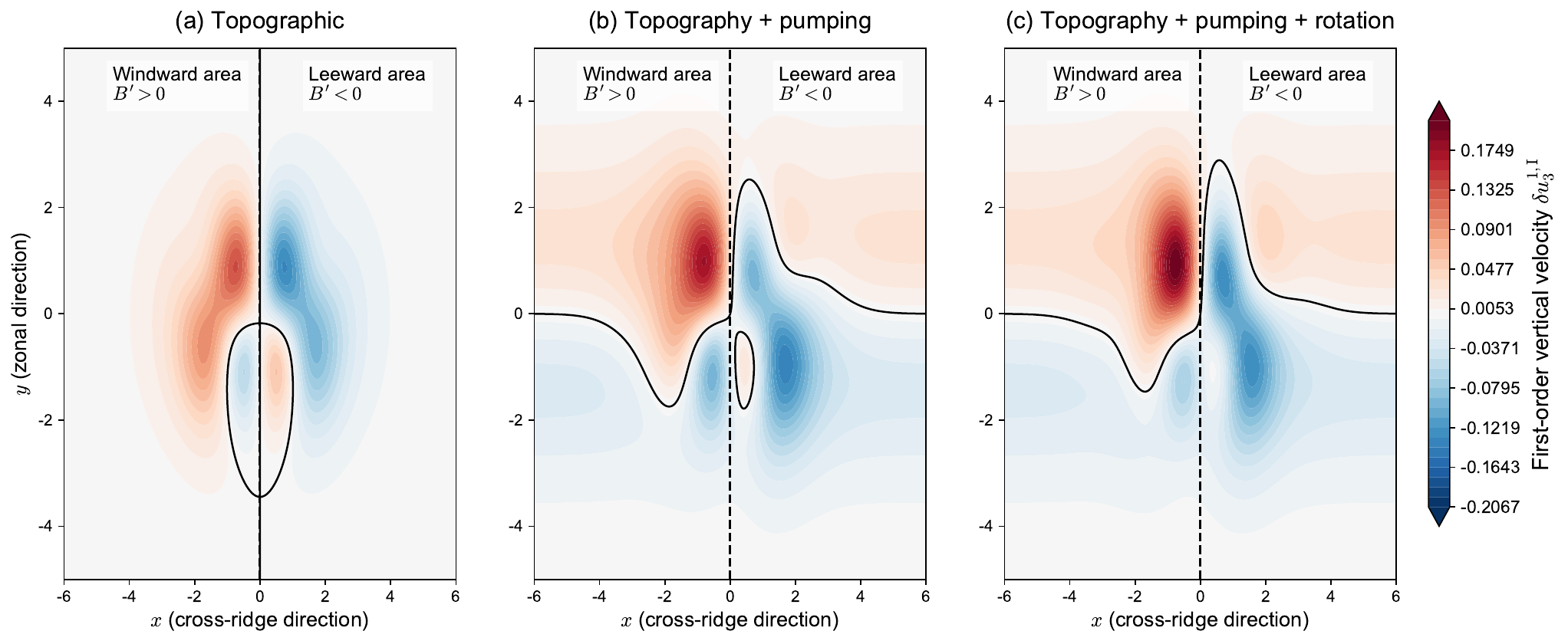}
    \caption{\textbf{Term-by-term superposition evolution of the first-order interior vertical velocity in the large-undulation model}}
    \label{fig:u3_components_full}
\end{figure}

To observe the superposition process of various mechanisms in the large-undulation model, Fig. \ref{fig:u3_components_full} decomposes the three components of the first-order interior vertical velocity $\delta u_3^{1,\mathrm{I}}$:

\begin{enumerate}
    \item \textbf{First-order topographic lifting (Fig. \ref{fig:u3_components_full}a):} 
This term reflects the geometric lifting of the fluid as it crosses the terrain. Overall, the fluid shows ascending on the windward slope ($x<0$) and descending on the leeward slope ($x>0$), with an antisymmetric distribution about the ridge ($x=0$). However, the extremum intensity in the upper half ($y>0$) is higher than that in the lower half ($y<0$), and closed zero-contour lines appear in the lower half. This phenomenon is attributed to the nonlinear shear of the first-order horizontal correction flow field ($\boldsymbol{u}_h^{1,\mathrm{I}}$) coupled with the background jet, leading to differences in the fluid on both sides.

    \item \textbf{Superimposed generalized Ekman suction (Fig. \ref{fig:u3_components_full}b):} 
After introducing the bottom friction induced by the metric tensor $\boldsymbol{H}_0$, the symmetry of the flow field is broken, the zero-contour lines become curved, and the suction intensity produces asymmetric offsets and magnitude increases on both sides. In addition, amplified by the thickness factor $\cos^{-1/2}\gamma$, the contour lines undergo transverse stretching at the steep parts of the ridge, reflecting the regulation of topographic geometry on the dissipation region.

    \item \textbf{Full effect superposition (Fig. \ref{fig:u3_components_full}c):}  
Finally, we introduce the rotational effect, which would not produce vertical velocity under flat terrain. However, comparing Fig. 2 and Fig. 3, we find that after adding the rotational term, the suction intensity is enhanced, and the black zero-contour line undergoes a significant transverse offset near the origin. This is mainly due to the coupling between the rotation operator $-\boldsymbol{E}_1$ and the thickness variation ($\cos^{-3/2}\gamma$).
\end{enumerate}

%
\appendix
\section{Appendices}

In this appendix, we present the proofs of Propositions \ref{D-prop2}–\ref{D-prop4} for the limiting system \eqref{D-3.45}. First, we give the equation satisfied by the macroscopic vorticity $\bar{\omega}=\nabla_h^\perp \cdot \bar{\boldsymbol{u}}_h$:
\begin{equation}\label{D-A1}
        \partial_t \bar{\omega} + (\bar{\boldsymbol u}_h \cdot \nabla_h) \bar{\omega}
        + \nabla_h^\perp \cdot \big(\sqrt{\tfrac{\nu}{2\cos\gamma}} (\boldsymbol{H}_0 -\cos^{-1}\gamma\boldsymbol{E}_1) \bar{\boldsymbol u}_h\big) + (\nabla_h B \cdot \nabla_h^\perp ) \,  D_t \bar{u}_3   = {0}, 
\end{equation}
To avoid the complexity caused by the geometric effect term $(\nabla_h B \cdot \nabla_h^\perp ) \,  D_t \bar{u}_3$, we rewrite the limiting system \eqref{D-3.45} into the following equivalent form:
\begin{equation}\label{D-A2}
    \partial_t (\boldsymbol{H}_0 \bar{\boldsymbol{u}}_h) 
    + \tilde{\omega} \bar{\boldsymbol{u}}_h^\perp + \nabla_h \Pi +\sqrt{\tfrac{\nu}{2\cos\gamma}} (\boldsymbol{H}_0 -\cos^{-1}\gamma\boldsymbol{E}_1) \bar{\boldsymbol u}_h = \boldsymbol{0},
\end{equation}
where the generalized vorticity $\tilde{\omega} = \nabla_h^\perp \cdot (\boldsymbol{H}_0 \bar{\boldsymbol{u}}_h)$ and $\Pi = \bar{p} + \frac{1}{2}|\bar{\boldsymbol{u}}_h|^2 + \frac{1}{2}\bar{u}_3^2$. 
Applying the 2D curl operator $\nabla_h^\perp \cdot$ to both sides of \eqref{D-A2} yields:
\begin{equation*}
    \partial_t \tilde{\omega} + (\bar{\boldsymbol{u}}_h \cdot \nabla_h)\tilde{\omega} + \nabla_h^\perp \cdot \big(\sqrt{\tfrac{\nu}{2\cos\gamma}} (\boldsymbol{H}_0 -\cos^{-1}\gamma\boldsymbol{E}_1) \bar{\boldsymbol u}_h\big) = 0.
\end{equation*}
Combining with
\[
	\nabla^\bot_h\gamma=-\tfrac{\cos^3\gamma}{\sin\gamma}\boldsymbol{E_1}\boldsymbol{H}\nabla_hB,
\]
we obtain the generalized vorticity equation:
\begin{align}\label{D-A3}
	&\partial_t\tilde{\omega}+(\bar {\boldsymbol u}_h\cdot\nabla_h) {\tilde{\omega}}
	+\sqrt{\tfrac{\nu}{2\cos\gamma}}\tilde{\omega}\\
	=&-\sqrt{\tfrac{\nu}{2\cos\gamma}}
	\cos^{-1}\gamma\nabla^\bot_hB\cdot
	\big(
	K_A\boldsymbol{E}
	-\tfrac{1}{2}\cos\gamma\boldsymbol{H}
	-\tfrac{3}{2}\cos^2\gamma\boldsymbol{E_1}\boldsymbol{H}
	\big)\bar{\boldsymbol{u}}_h\notag.
\end{align}

\subsection{Proof of Proposition \ref{D-prop2}}

First, we calculate the $L^2$ energy estimate for the velocity field $\bar{\boldsymbol{u}}_h$. Taking the $L^2$ inner product of both sides of equation \eqref{D-3.45} over $\mathbb{R}^2$ yields:
\begin{multline*}
    \frac{1}{2}\frac{d}{dt} \|\bar{\boldsymbol{u}}_h\|_{L^2(\mathbb{R}^2)}^2 + \sqrt{\tfrac{\nu}{2}} \int_{\mathbb{R}^2} \cos^{-\frac12}\gamma \big(|\bar{\boldsymbol{u}}_h|^2+(\nabla_h B \cdot \bar{\boldsymbol{u}}_h)^2 \big) \,dxdy \\+ \int_{\mathbb{R}^2} \bar{\boldsymbol{u}}_h \cdot \nabla_h B  \big( \partial_t \bar{u}_3 + (\bar{\boldsymbol{u}}_h \cdot \nabla_h)\bar{u}_3 \big) \,dxdy = 0.
\end{multline*}
Due to the incompressibility condition of $\bar{\boldsymbol{u}}_h$ and $\bar{u}_3 = \nabla_h B \cdot \bar{\boldsymbol{u}}_h$, the topographic feedback term simplifies to:
\begin{equation*}
    \int_{\mathbb{R}^2} \bar{\boldsymbol{u}}_h \cdot \nabla_h B \big( \partial_t \bar{u}_3 + (\bar{\boldsymbol{u}}_h \cdot \nabla_h)\bar{u}_3 \big) \,dxdy = \frac{1}{2} \frac{d}{dt} \|\bar{u}_3\|_{L^2(\mathbb{R}^2)}^2.
\end{equation*}
Substituting this back into the energy equality and combining with $\cos^{-\frac12}\gamma> 1$, we obtain:
\begin{equation}\label{D-A4}
    \frac{1}{2}\frac{d}{dt}  \|\bar{\boldsymbol{u}}\|_{L^2(\mathbb{R}^2)}^2   + \sqrt{\tfrac{\nu}{2}}  \|\bar{\boldsymbol{u}}\|_{L^2(\mathbb{R}^2)}^2 \leqslant 0.
\end{equation}
Finally, we get the decay estimate for $\bar{\boldsymbol{u}}$:
\begin{equation}\label{D-A5}
    \|\bar{\boldsymbol{u}}(t)\|_{L^2(\mathbb{R}^2)}^2  \leqslant  \|\bar{\boldsymbol{u}}_0\|_{L^2(\mathbb{R}^2)}^2  \mathrm{e}^{-{\sqrt{2\nu}}t}.
\end{equation}

Next, we perform vorticity estimation for the flow field. First, we analyze the generalized vorticity $\tilde{\omega}$.
Multiplying both sides of equation \eqref{D-A3} by $\tilde{\omega}$ and integrating yields:
	\begin{align}\label{D-A6}
		&	\frac12\frac{d}{dt}\lVert\tilde{\omega}\rVert_{L^2(\mathbb{R}^2)}^2+\sqrt{\tfrac{\nu}{2}}\lVert\tilde{\omega}\rVert_{L^2(\mathbb{R}^2)}^2
		\\\notag
		&\leqslant\sqrt{\tfrac{\nu}{2}}\Big|\int_{\mathbb{R}^2}\cos^{-1}\gamma\nabla^\bot_hB\cdot
        \big(
		K_A\boldsymbol{E}
		-\tfrac{1}{2}\cos\gamma\boldsymbol{H}
		-\tfrac{3}{2}\cos^{2}\gamma\boldsymbol{E_1}\boldsymbol{H}
		\big)\bar{\boldsymbol{u}}_h \cdot \tilde{\omega}  \,dx{d}y\Big|
		\notag.
	\end{align}
For the low-order linear coupling term induced by non-flat geometry that appears in the error estimate, we denote:
\begin{equation*}
    \boldsymbol{A} = K_A\boldsymbol{E} -\tfrac{1}{2}\cos\gamma\boldsymbol{H} -\tfrac{3}{2}\cos^{2}\gamma\boldsymbol{E_1}\boldsymbol{H}.
\end{equation*}
Using the triangle inequality for operator norms, and noting that $\|\boldsymbol{E}\| = \|\boldsymbol{E_1}\| = 1$ and $\cos\gamma \leqslant 1$, we can bound the spectral norm of $\boldsymbol{A}$ as:
\begin{align}\label{D-A7}
    \|\boldsymbol{A}\| 
    &\leqslant|K_A| + \tfrac{1}{2}\cos\gamma\|\boldsymbol{H}\| + \tfrac{3}{2}\cos^2\gamma\|\boldsymbol{H}\| \\
    &\leqslant |K_A| + 2\max|\lambda_{\pm}| \cos\gamma ,\notag
\end{align}
where $\lambda_{\pm}$ are the eigenvalues of the matrix $\boldsymbol{H}$:
\begin{equation}\label{D-A8}
\lambda_{\pm} = \tfrac{1}{2}\Big(\Delta_h B \pm \sqrt{(\Delta_h B)^2 - 4\cos^{-4}\gamma K_G}\Big).
\end{equation}
Although we cannot write the horizontal Laplacian $\Delta_h B$ as an explicit equality in terms of mean curvature and Gaussian curvature when dealing with a priori estimates involving the second derivative of topography $\Delta_h B$, we can give its upper and lower bounds using the principal curvatures of the boundary geometry, as stated in the following lemma:
\begin{lemma}\label{D-lem1}
Under the geometric assumptions of Proposition \ref{D-prop2},
if we let $k_n$ be the normal curvature along the topographic gradient direction $\boldsymbol{v}_B = \nabla_h B / |\nabla_h B|$, then the horizontal Laplacian operator $\Delta_h B$ of the topography has the following geometric decomposition:
\begin{equation}\label{D-A9}
    \Delta_h B = k_n (\cos^{-3}\gamma - \cos^{-1}\gamma) + 2K_A \cos^{-1}\gamma,
\end{equation}
and satisfies:
\begin{equation}\label{D-A10}
    |\Delta_h B| \leqslant 2C_0  \cos^{-3}\gamma  .
\end{equation}
\end{lemma}
\begin{proof}
First, the mean curvature $K_A$ of the boundary surface $z=B(x,y)$ has the following representation:
\begin{equation}\label{D-A11}
    2K_A = \nabla_h \cdot \big( \cos\gamma \nabla_h B \big)= \cos\gamma \Delta_h B -\cos^3\gamma \nabla_h B^T\boldsymbol{H}\nabla_h B.
\end{equation}
The normal curvature $k_n$ along the topographic gradient direction $\boldsymbol{v}_B = \nabla_h B / |\nabla_h B|$ can be expressed in terms of its first fundamental form $\mathrm{I}$ and second fundamental form $\mathrm{II}$ as:
\begin{equation}\label{D-A12}
    k_n = \frac{\mathrm{II}(\boldsymbol{v}_B,\boldsymbol{v}_B)}{\mathrm{I}(\boldsymbol{v}_B,\boldsymbol{v}_B)} = \frac{\cos\gamma \, \boldsymbol{v}_B^T \boldsymbol{H} \boldsymbol{v}_B}{1 + (\boldsymbol{v}_B \cdot \nabla_h B)^2} = \cos^3\gamma \frac{\nabla_h B^T \boldsymbol{H} \nabla_h B}{\tan^2\gamma}.
\end{equation}
Then, from \eqref{D-A11}–\eqref{D-A12}, we can invert to obtain the form of $\Delta_h B$:
\begin{equation*}
    \Delta_h B = k_n (\cos^{-3}\gamma - \cos^{-1}\gamma) + 2K_A \cos^{-1}\gamma.
\end{equation*} 
It is known that the normal curvature lies between the two principal curvatures (i.e., $k_2 \leqslant k_n \leqslant k_1$).
Since $\cos^{-3}\gamma \geqslant \cos^{-1}\gamma > 1$ and using the triangle inequality, we get:
\begin{equation*}
    |\Delta_h B| 
    =|k_n(\cos^{-3}\gamma - \cos^{-1}\gamma) + (k_1+k_2) \cos^{-1}\gamma|
    \leqslant 2 C_0 \cos^{-3}\gamma.
\end{equation*} 
\end{proof}
After obtaining the bound for the horizontal Laplacian $\Delta_h B$, using the basic inequality $\sqrt{x^2+y^2} \leqslant |x|+|y|$, we naturally derive the following estimate for the eigenvalues of matrix $\boldsymbol{H}$:
\begin{align}\label{D-A13}
|\lambda_{\pm}| &=\tfrac{1}{2}\Big|\Delta_h B \pm \sqrt{(\Delta_h B)^2 - 4\cos^{-4}\gamma K_G}\Big|\\
&\leqslant \tfrac{1}{2}|\Delta_h B| + \tfrac{1}{2}\Big( |\Delta_h B| + 2\cos^{-2}\gamma \sqrt{|K_{\mathrm{G}}|} \Big)\notag \\
    &= |\Delta_h B| + \cos^{-2}\gamma \sqrt{|K_{\mathrm{G}}|}.\notag
\end{align}
Note that the Gaussian curvature satisfies $\sqrt{|K_{\mathrm{G}}|} = \sqrt{|k_1 k_2|} \leqslant  C_0 $, and the mean curvature satisfies ${|K_{A}|} \leqslant \frac12 C_0$. 
Thus, we obtain the bound for the spectral norm of $\boldsymbol{A}$:
\begin{equation}\label{D-A14}
\|\boldsymbol{A}\| \leqslant \tfrac12 C_0 \big(1+ 12  \cos^{-2}\gamma  \big).
\end{equation}

Next, combining equation \eqref{D-A14} and the geometric constraint \eqref{D-3.46}, we further estimate inequality \eqref{D-A6}:
	\begin{align}\label{D-A15}
		&	\frac{1}{2}\frac{d}{dt}\lVert\tilde{\omega}\rVert_{L^2(\mathbb{R}^2)}^2+\sqrt{\tfrac{\nu}{2}}\lVert\tilde{\omega}\rVert_{L^2(\mathbb{R}^2)}^2
		\\\notag
		&\leqslant \tfrac{C_0}{2}  \sqrt{\tfrac{\nu}{2}}
        \int_{\mathbb{R}^2}|\nabla_h B|
        \, \big(1+12  \cos^{-2}\gamma  \big)
        |\bar{\boldsymbol{u}}_h| |\tilde{\omega}| \,dx{d}y
		\notag\\
        &\lesssim
        \sqrt{\tfrac{\nu}{2}}\,
        \int_{\mathbb{R}^2}
        |\bar{\boldsymbol{u}}_h| |\tilde{\omega}| \,dx{d}y
		\notag\\
        &\leqslant
        \tfrac{1}{2} \sqrt{\tfrac{\nu}{2}} \|\tilde{\omega}\|_{L^2(\mathbb{R}^2)}^2 + \tfrac{1}{2} \sqrt{\tfrac{\nu}{2}}\,  \|\bar{\boldsymbol{u}}_h\|_{L^2(\mathbb{R}^2)}^2.\notag
	\end{align}

Furthermore, combining the decay estimate for the horizontal velocity field $\|\bar{\boldsymbol{u}}_h\|_{L^2}^2$ \eqref{D-A5} and the definition of generalized vorticity $\tilde{\omega} = \nabla_h^\perp \cdot (\boldsymbol{H}_0 \bar{\boldsymbol{u}}_h)$, we easily derive the global time decay estimate for the generalized vorticity:
\begin{align*}
    \lVert\tilde{\omega}(t)\rVert_{L^2(\mathbb{R}^2)}^2 
    \lesssim&  \Big( \lVert\tilde{\omega}_0\rVert_{L^2(\mathbb{R}^2)}^2 +   \|\bar{\boldsymbol{u}}_0\|_{L^2(\mathbb{R}^2)}^2 \Big) \mathrm{e}^{-\sqrt{\frac{\nu}{2}}t}\\
    \lesssim&  \Big(\|\bar{\omega}_0\|^2_{L^2(\mathbb{R}^2)}+\|\bar{\boldsymbol{u}}_{0,h}\|^2_{L^2(\mathbb{R}^2)}   \Big) \mathrm{e}^{-\sqrt{\frac{\nu}{2}}t}.
\end{align*} 

Finally, we need to reduce the exponential decay result of the generalized vorticity $\tilde{\omega}$ to the vorticity $\bar{\omega} = \nabla_h^\perp \cdot \bar{\boldsymbol{u}}_h$. According to the definition of generalized vorticity, we have:
\begin{equation*}
    \|\bar{\omega}(t)\|^2_{L^2(\mathbb{R}^2)} \lesssim \|\tilde{\omega}(t)\|^2_{L^2(\mathbb{R}^2)} + \|\bar{\boldsymbol{u}}_h(t)\|^2_{L^2(\mathbb{R}^2)}\lesssim \Big(\|\bar{\omega}_0\|^2_{L^2(\mathbb{R}^2)}+\|\bar{\boldsymbol{u}}_{0,h}\|^2_{L^2(\mathbb{R}^2)}   \Big) \mathrm{e}^{-\sqrt{\frac{\nu}{2}}t}.
\end{equation*}

\subsection{Proof of Proposition \ref{D-prop3}}
    First, we establish estimates for the generalized vorticity in $L^s$ spaces ($2 \le s < \infty$). Multiplying both sides of the generalized vorticity equation \eqref{D-A3} by $\tilde{\omega}|\tilde{\omega}|^{s-2}$ and integrating over the whole space $\mathbb{R}^2$. Since the velocity field satisfies $\nabla_h \cdot \bar{\boldsymbol{u}}_h = 0$ and the boundary $B(x,y)$ is smooth, the $L^s$ norm of the generalized vorticity satisfies the following inequality:
    \begin{equation}\label{D-A16}
        \frac{1}{s}\frac{d}{dt} \|\tilde{\omega}\|^s_{L^s(\mathbb{R}^2)} + \sqrt{\frac{\nu}{2}} \|\tilde{\omega}\|^s_{L^s(\mathbb{R}^2)} \lesssim \|\tilde{\omega}\|^{s-1}_{L^s(\mathbb{R}^2)} \|\bar{\boldsymbol{u}}_h\|_{L^s(\mathbb{R}^2)}.
    \end{equation}
    Canceling the common factor $\|\tilde{\omega}\|^{s-1}_{L^s(\mathbb{R}^2)}$ and applying Gronwall's lemma, we solve:
    \begin{equation}\label{D-A17}
        \|\tilde{\omega}(t)\|_{L^s(\mathbb{R}^2)} \lesssim \mathrm{e}^{-\sqrt{\frac{\nu}{2}} t} \|\tilde{\omega}_0\|_{L^s(\mathbb{R}^2)} + \int_0^t \mathrm{e}^{-\sqrt{\frac{\nu}{2}} (t-\tau)} \|\bar{\boldsymbol{u}}_h(\tau)\|_{L^s(\mathbb{R}^2)} \, d\tau.
    \end{equation}

    Specifically choosing $s=4$, combining with Proposition \ref{D-prop2}, and using the 2D Gagliardo-Nirenberg interpolation inequality, we get:
    \begin{equation}\label{D-A18}
        \|\bar{\boldsymbol{u}}_h(t)\|_{L^4(\mathbb{R}^2)} \lesssim \|\bar{\boldsymbol{u}}_h\|_{L^2(\mathbb{R}^2)}^{1/2} \|\nabla_h \bar{\boldsymbol{u}}_h\|_{L^2(\mathbb{R}^2)}^{1/2} \lesssim \mathcal{I}_0 \, \mathrm{e}^{-\frac{1}{2}\sqrt{\frac{\nu}{2}}t}.
    \end{equation}
    Substituting \eqref{D-A18} back into the integral \eqref{D-A17}, and using the initial value $\|\tilde{\omega}_0\|_{L^4(\mathbb{R}^2)} \lesssim \|\tilde{\omega}_0\|_{L^2(\mathbb{R}^2)}^{1/2}\|\tilde{\omega}_0\|_{L^\infty(\mathbb{R}^2)}^{1/2} \lesssim \mathcal{I}_0$, we infer that the $L^4$ norm of $\tilde{\omega}$ satisfies:
    \begin{equation}\label{D-A19}
        \|\tilde{\omega}(t)\|_{L^4(\mathbb{R}^2)} \lesssim \mathcal{I}_0 \, \mathrm{e}^{-\frac{1}{2}\sqrt{\frac{\nu}{2}} t}.
    \end{equation}

    Subsequently, according to the 2D Sobolev embedding $W^{1,4}(\mathbb{R}^2) \hookrightarrow L^\infty(\mathbb{R}^2)$, we obtain:
\begin{equation}\label{D-A20}
    \|\bar{\boldsymbol{u}}_h(t)\|_{L^\infty(\mathbb{R}^2)} \lesssim \|\bar{\boldsymbol{u}}_h(t)\|_{W^{1,4}(\mathbb{R}^2)} \lesssim \|\tilde{\omega}(t)\|_{L^4(\mathbb{R}^2)} + \|\bar{\boldsymbol{u}}_h(t)\|_{L^2(\mathbb{R}^2)} \lesssim \mathcal{I}_0 \, \mathrm{e}^{-\frac{1}{2}\sqrt{\frac{\nu}{2}} t}.
\end{equation}

Finally, integrating equation \eqref{D-A3} along the characteristic lines of the macroscopic flow field and taking the $L^\infty$ norm, the geometric perturbation induced by topography on the right-hand side is controlled by the $L^\infty$ norm of the velocity field, so we have:
    \begin{equation*}
        \|\tilde{\omega}(t)\|_{L^\infty(\mathbb{R}^2)} \lesssim \mathrm{e}^{-\sqrt{\frac{\nu}{2}} t} \|\tilde{\omega}_0\|_{L^\infty(\mathbb{R}^2)} + \int_0^t \mathrm{e}^{-\sqrt{\frac{\nu}{2}} (t-\tau)} \|\bar{\boldsymbol{u}}_h(\tau)\|_{L^\infty(\mathbb{R}^2)} \, d\tau.
    \end{equation*}
    Substituting equation \eqref{D-A20} into the integral term above, we obtain the exponential decay estimate for the generalized vorticity in the $L^\infty$ framework:
    \begin{equation*}
        \|\tilde{\omega}(t)\|_{L^\infty(\mathbb{R}^2)} \lesssim   \mathcal{I}_0 \, \mathrm{e}^{-\frac{1}{2}\sqrt{\frac{\nu}{2}} t}.
    \end{equation*}

\subsection{Proof of Proposition \ref{D-prop4}}
    According to the definition of generalized vorticity, the high-order Sobolev norms of the velocity field are controlled by the generalized vorticity. In particular, for $H^2$ and $H^3$ spaces, we have $\|\bar{\boldsymbol{u}}_h\|_{H^2} \lesssim \|\tilde{\omega}\|_{H^1} + \|\bar{\boldsymbol{u}}_h\|_{L^2}$ and $\|\bar{\boldsymbol{u}}_h\|_{H^3} \lesssim \|\tilde{\omega}\|_{H^2} + \|\bar{\boldsymbol{u}}_h\|_{L^2}$. 
    
    Next, applying the $\nabla_h^2$ operator to the generalized vorticity equation and taking the $L^2$ inner product with $\nabla_h^2 \tilde{\omega}$ yields:
\begin{align}\label{D-A21}
    &\frac{1}{2} \frac{d}{dt} \|\tilde{\omega}\|_{H^2(\mathbb{R}^2)}^2 + \sqrt{\frac{\nu}{2}} \|\tilde{\omega}\|_{H^2(\mathbb{R}^2)}^2  \\
    &\lesssim \left| \int_{\mathbb{R}^2}[\nabla_h^2, \bar{\boldsymbol{u}}_h \cdot \nabla_h] \tilde{\omega} \cdot \nabla_h^2 \tilde{\omega} \,dxdy \right| +  \|\bar{\boldsymbol{u}}_h\|_{H^2(\mathbb{R}^2)} \|\tilde{\omega}\|_{H^2(\mathbb{R}^2)} \notag \\
    &\lesssim \|\nabla_h \bar{\boldsymbol{u}}_h\|_{L^\infty(\mathbb{R}^2)} \|\tilde{\omega}\|_{H^2(\mathbb{R}^2)}^2 +\big( \|\tilde{\omega}\|_{H^1(\mathbb{R}^2)} + \|\bar{\boldsymbol{u}}_h\|_{L^2(\mathbb{R}^2)} \big) \|\tilde{\omega}\|_{H^2(\mathbb{R}^2)}\notag\\
    &\lesssim \|\nabla_h \bar{\boldsymbol{u}}_h\|_{L^\infty(\mathbb{R}^2)} \|\tilde{\omega}\|_{H^2(\mathbb{R}^2)}^2 +
    \|\tilde{\omega}\|_{L^2(\mathbb{R}^2)}^{1/2} \|\tilde{\omega}\|_{H^2(\mathbb{R}^2)}^{3/2} +  \|\bar{\boldsymbol{u}}_h\|_{L^2(\mathbb{R}^2)} \|\tilde{\omega}\|_{H^2(\mathbb{R}^2)}
    \notag\\
    &\lesssim \frac{1}{2}\sqrt{\frac{\nu}{2}} \|\tilde{\omega}\|_{H^2(\mathbb{R}^2)}^2+\|\nabla_h \bar{\boldsymbol{u}}_h\|_{L^\infty(\mathbb{R}^2)} \|\tilde{\omega}\|_{H^2(\mathbb{R}^2)}^2 
    +   \|\tilde{\omega}\|_{L^2(\mathbb{R}^2)}^2 + \|\bar{\boldsymbol{u}}_h\|_{L^2(\mathbb{R}^2)}^2 
    ,\notag
\end{align}
where $[\cdot, \cdot]$ denotes the commutator. 

Combining the results of Proposition \ref{D-prop2} and introducing the functional $Y(t) \triangleq \mathrm{e} + \|\tilde{\omega}(t)\|_{H^2(\mathbb{R}^2)}^2$, we can scale \eqref{D-A21} into a differential inequality for $Y(t)$:
\begin{equation}\label{D-A22}
    \frac{d}{dt} Y(t) + \sqrt{\frac{\nu}{2}} (Y(t) - \mathrm{e}) \lesssim \|\nabla_h \bar{\boldsymbol{u}}_h\|_{L^\infty(\mathbb{R}^2)} Y(t) + \mathrm{e}^{-\sqrt{\frac{\nu}{2}} t}.
\end{equation}
Since $Y(t) \ge \mathrm{e} > 1$, we can further scale the above inequality as:
\begin{equation}\label{D-A23}
    \frac{d}{dt} Y(t) \lesssim \big( \|\nabla_h \bar{\boldsymbol{u}}_h\|_{L^\infty(\mathbb{R}^2)} + \mathrm{e}^{-\sqrt{\frac{\nu}{2}} t} \big) Y(t).
\end{equation}
To compute the $L^\infty$ estimate for $\nabla_h \bar{\boldsymbol{u}}_h$, we use the proof idea from Ref.\refcite{BKM1984} to derive a logarithmic Sobolev inequality. 

\begin{lemma}\label{D-lem2}
Assume that the 2D divergence-free velocity field is $\bar{\boldsymbol{u}}_h = \nabla_h^\perp \psi \in H^3(\mathbb{R}^2)$ and the generalized vorticity is $\tilde{\omega} = \nabla_h^\perp \cdot (\boldsymbol{H}_0 \bar{\boldsymbol{u}}_h)\in L^\infty(\mathbb{R}^2)$. 
If the supremum of the principal curvatures of the topography over the whole space is $C_0 <\infty$, then the following logarithmic Sobolev inequality holds:
\begin{equation*}
    \|\nabla_h \bar{\boldsymbol{u}}_h\|_{L^\infty(\mathbb{R}^2)} \lesssim (1 + C_0) \Big( \|\bar{\boldsymbol{u}}_h\|_{L^2(\mathbb{R}^2)} + \|\tilde{\omega}\|_{L^\infty(\mathbb{R}^2)} \ln \big(\mathrm{e} + \|\bar{\boldsymbol{u}}_h\|_{H^3(\mathbb{R}^2)}\big) \Big).
\end{equation*}
\end{lemma}

\begin{proof}
Since the macroscopic flow field satisfies $\nabla_h \cdot \bar{\boldsymbol{u}}_h = 0$, there exists a stream function $\psi$ such that $\bar{\boldsymbol{u}}_h = \nabla_h^\perp \psi$. Substituting this into the definition of generalized vorticity, we obtain a variable-coefficient elliptic Poisson equation for $\psi$:
\begin{equation*}
    \mathcal{L} \psi := \nabla_h^\perp \cdot (\boldsymbol{H}_0 \nabla_h^\perp \psi) = \tilde{\omega}.
\end{equation*}

Using the inverse operator $\mathcal{L}^{-1}$ of the variable-coefficient operator $\mathcal{L}$, we can explicitly express the horizontal velocity gradient as:
\begin{equation*}
    \nabla_h \bar{\boldsymbol{u}}_h = \big( \nabla_h \nabla_h^\perp \mathcal{L}^{-1} \big) \tilde{\omega} := \mathcal{T} \tilde{\omega}.
\end{equation*}
Here $\mathcal{T}$ is a zero-order variable-coefficient differential operator. To handle this variable-coefficient singular integral operator $\mathcal{T} = \nabla_h \nabla_h^\perp \mathcal{L}^{-1}$, we introduce the Littlewood-Paley dyadic decomposition. Let $\{\Delta_q\}_{q \geqslant -1}$ be a family of smooth cutoff operators in the frequency domain. We decompose the full-space $L^\infty$ norm of the velocity gradient into three frequency bands: low frequency ($q \leqslant 0$), middle frequency ($0 < q \leqslant N$), and high frequency ($q > N$), where $N \in \mathbb{N}^+$ is to be determined. Then we have:
\begin{align*}
    \|\nabla_h \bar{\boldsymbol{u}}_h\|_{L^\infty}  &\leqslant\sum_{q=-1}^0 \|\Delta_q \nabla_h \bar{\boldsymbol{u}}_h\|_{L^\infty} + \sum_{q=1}^N \|\Delta_q \nabla_h \bar{\boldsymbol{u}}_h\|_{L^\infty} + \sum_{q=N+1}^\infty \|\Delta_q \nabla_h \bar{\boldsymbol{u}}_h\|_{L^\infty} \\
    &:= I_{low} + I_{mid} + I_{high}.
\end{align*}

\textbf{Step 1: Low-frequency estimate ($I_{low}$)} \\
For the low-frequency part ($q \leqslant 0$), according to Bernstein's inequality on the 2D full space $\mathbb{R}^2$, we get:
\begin{equation*}
    I_{low} \leqslant \sum_{q=-1}^0 2^q \|\Delta_q \bar{\boldsymbol{u}}_h\|_{L^\infty} \leqslant \sum_{q=-1}^0 2^{2q} \|\Delta_q \bar{\boldsymbol{u}}_h\|_{L^2} \leqslant C_1 \|\bar{\boldsymbol{u}}_h\|_{L^2}.
\end{equation*}

\textbf{Step 2: Middle-frequency estimate ($I_{mid}$) and curvature control} \\
For the middle-frequency part ($0 < q \leqslant N$), we use the relation with the generalized vorticity $\nabla_h \bar{\boldsymbol{u}}_h = \mathcal{T} \tilde{\omega}$. Since the principal curvatures of the bottom topography are bounded, the metric tensor satisfies $\|\nabla \boldsymbol{H}_0\|_{L^\infty} \lesssim C_0$. According to the commutator theory of singular integrals, the operator $\mathcal{T}$ is uniformly $L^\infty$-bounded within each individual dyadic frequency block, and its operator norm is dominated by the Lipschitz constant. Therefore:
\begin{equation*}
    I_{mid} = \sum_{q=1}^N \|\Delta_q (\mathcal{T} \tilde{\omega})\|_{L^\infty} \leqslant \sum_{q=1}^N C(1+C_0) \|\Delta_q \tilde{\omega}\|_{L^\infty} \leqslant C_2 (1+C_0) N \|\tilde{\omega}\|_{L^\infty}.
\end{equation*}

\textbf{Step 3: High-frequency estimate ($I_{high}$)} \\
For the high-frequency part ($q > N$), using Bernstein's inequality again and the high-order regularity of the flow field, we get:
\begin{align*}
    I_{high} \leqslant &\sum_{q > N} 2^{2q} \|\Delta_q \bar{\boldsymbol{u}}_h\|_{L^2} = \sum_{q > N} 2^{-q} \big(2^{3q} \|\Delta_q \bar{\boldsymbol{u}}_h\|_{L^2}\big) \\
    \leqslant& \sum_{q > N} 2^{-q} \|\bar{\boldsymbol{u}}_h\|_{H^3} \leqslant C_3 2^{-N} \|\bar{\boldsymbol{u}}_h\|_{H^3}.
\end{align*}

Combining the above three steps, we obtain an inequality in terms of the truncation order $N$:
\begin{equation*}
    \|\nabla_h \bar{\boldsymbol{u}}_h\|_{L^\infty} \leqslant C \Big( \|\bar{\boldsymbol{u}}_h\|_{L^2} + (1+C_0) N \|\tilde{\omega}\|_{L^\infty} + 2^{-N} \|\bar{\boldsymbol{u}}_h\|_{H^3} \Big).
\end{equation*}
Without loss of generality, we take $N \approx \log_2 \big( \mathrm{e} + \|\bar{\boldsymbol{u}}_h\|_{H^3} \big)$. Finally, we obtain the logarithmic generalized BKM inequality. 
\end{proof}

Using the result of Lemma \ref{D-lem2} and the definition of $Y(t)$, we get:
\begin{align}\label{D-A24}
    \|\nabla_h \bar{\boldsymbol{u}}_h(t)\|_{L^\infty(\mathbb{R}^2)} 
    \lesssim  (1+C_0)\big( 1 + \ln Y(t) \big) \mathrm{e}^{-\frac{1}{2}\sqrt{\frac{\nu}{2}} t}.
\end{align}
Substituting \eqref{D-A24} back into \eqref{D-A23} and dividing both sides by $Y(t)$ yields:
\begin{equation}\label{D-A25}
    \frac{d}{dt} \big( \ln Y(t) \big) \lesssim (1+C_0)\big( 1 + \ln Y(t) \big) \mathrm{e}^{-\frac{1}{2}\sqrt{\frac{\nu}{2}} t} .
\end{equation}
Note that the coefficient time function on the right-hand side of the above differential inequality is integrable over $[0, \infty)$. Applying Gronwall's lemma to equation \eqref{D-A25} shows that $\ln Y(t)$ is uniformly bounded, i.e., there exists a constant $C > 0$ such that:
\begin{equation}\label{D-A26}
\ln Y(t) \le C < \infty, \quad \forall t \ge 0.
\end{equation}
This implies:
\begin{equation}\label{D-A27}
    \|\nabla_h \bar{\boldsymbol{u}}_h(t)\|_{L^\infty(\mathbb{R}^2)} \lesssim \big( 1 + C \big)\mathrm{e}^{-\frac{1}{2}\sqrt{\frac{\nu}{2}} t} .
\end{equation}
as well as that $\|\tilde{\omega}\|_{H^2}$ and $\|\bar{\boldsymbol{u}}_h\|_{H^3}$ are uniformly bounded. Based on this fact, and retaining the damping term on the left-hand side of inequality \eqref{D-A21}, we apply Gronwall's lemma again to obtain:
\begin{equation*}
    \|\tilde{\omega}(t)\|_{H^2(\mathbb{R}^2)} \lesssim \mathrm{e}^{-\frac{1}{2}\sqrt{\frac{\nu}{2}} t},
\end{equation*}
and consequently:
\begin{equation*}
    \|\bar{\boldsymbol{u}}_h\|_{H^3(\mathbb{R}^2)} \lesssim \mathrm{e}^{-\frac{1}{2}\sqrt{\frac{\nu}{2}} t}.
\end{equation*}

%
%
%

\section*{Acknowledgment}
This work was partially supported by the National Key R\&D Program of China (No. 2020YFA072500), Central Guidance for Local Science and Technology Development Fund (No. ZYYD2026ZY11), the Xinjiang Talent Development Fund (No. XJRC-2025-KJ-PY-KJLJ-105), and the Innovation Project of Excellent Doctoral Students of Xinjiang University (No. XJU2024BS038)
The authors thank the anonymous reviewers for their helpful comments.

\end{document}